\documentclass[10pt]{article}
\usepackage{amsmath,amssymb,amsthm}
\begin{document}
\theoremstyle{plain}
\newtheorem{pro1}{Proposition}[section]
\newtheorem{cor1}[pro1]{Corollary}
\newtheorem{lem1}[pro1]{Lemma}
\newtheorem{thm1}[pro1]{Theorem}
\theoremstyle{definition}
\newtheorem{def1}[pro1]{Definition}
\newtheorem{fac1}[pro1]{Fact}
\newtheorem{rem1}[pro1]{Remark}
\newtheorem{ex1}[pro1]{Example}
\renewcommand{\labelenumi}{(\roman{enumi})}
\font\germ=eufm10
\def\al{\alpha}
\def\beneme{\begin{enumerate}}
\def\beq{\begin{equation}}
\def\beqn{\begin{eqnarray}}
\def\beqnn{\begin{eqnarray*}}
\def\bigsl{{\hbox{\fontD \char'54}}}
\def\cd{\cdots}
\def\del{\delta}
\def\Del{\Delta}
\def\ei{e_i}
\def\eit{\tilde{e}_i}
\def\eneme{\end{enumerate}}
\def\ep{\epsilon}
\def\eeq{\end{equation}}
\def\eeqn{\end{eqnarray}}
\def\eeqnn{\end{eqnarray*}}
\def\fit{\tilde{f}_i}
\def\ft{\tilde{f}}
\def\ge{{\mathfrak g}}
\def\gl{\hbox{\germ gl}}
\def\hom{{\hbox{Hom}}}
\def\ify{\infty}
\def\io{\iota}
\def\kp{k^{(+)}}
\def\km{k^{(-)}}
\def\llra{\relbar\joinrel\relbar\joinrel\relbar\joinrel\rightarrow}
\def\lan{\langle}
\def\lar{\longrightarrow}
\def\lm{\lambda}
\def\Lm{\Lambda}
\def\mapright#1{\smash{\mathop{\longrightarrow}\limits^{#1}}}
\def\nd{\noindent}
\def\nn{\nonumber}
\def\ot{\otimes}
\def\op{\oplus}
\def\opi{\ovl\pi_{\lm}}
\def\ovl{\overline}
\def\plm{\Psi^{(\lm)}_{\io}}
\def\qq{\qquad}
\def\q{\quad}
\def\qed{\hfill\framebox[3mm]{}}
\def\QQ{\mathbb Q}
\def\qi{q_i}
\def\qii{q_i^{-1}}
\def\ran{\rangle}
\def\rlm{r_{\lm}}
\def\ssl{\hbox{\germ sl}}
\def\slh{\widehat{\ssl_2}}
\def\ti{t_i}
\def\tii{t_i^{-1}}
\def\til{\tilde}
\def\tt{{\hbox{\germ{t}}}}
\def\ttt{\hbox{\germ t}}
\def\uq{U_q(\ge)}
\def\uqm{U^-_q(\ge)}
\def\uqp{U^+_q(\ge)}
\def\uqmq{{U^-_q(\ge)}_{\bf Q}}
\def\uqpm{U^{\pm}_q(\ge)}
\def\uqmp{U^{\mp}_q(\ge)}
\def\uqq{U_{\bf Q}^-(\ge)}
\def\uqz{U^-_{\bf Z}(\ge)}
\def\util{\tilde{U}_q(\ge)}
\def\vep{\varepsilon}
\def\vp{\varphi}
\def\vpi{\varphi^{-1}}
\def\xii{\xi^{(i)}}
\def\Xiioi{\Xi_{\io}^{(i)}}
\def\wtil{\widetilde}
\def\what{\widehat}
\def\wpi{\widehat\pi_{\lm}}
\def\ZZ{\mathbb Z}
%\def\zlm{\ZZ^{\ify}_{\io}[\lm]}
%\tableofcontents
\font\germ=eufm10
\def\ssl{\hbox{\germ sl}}
\def\slh{\widehat{\ssl_2}}
\title{\large \bf Polyhedral Realizations of Crystal Bases for 
\\Modified Quantum Algebras of Arbitrary Rank 2 Cases}
\author{HOSHINO Ayumu\thanks{e-mail address:
a-hoshin@mm.sophia.ac.jp} \\
\\
\normalsize Department of Mathematics, Sophia University, \\
\normalsize Tokyo 102-8554, JAPAN}
\date{}
\maketitle
\begin{abstract}
We describe the crystal bases of the modified quantum algebras
and give the explicit form of the highest (or lowest) weight vector of its
connected component $B_0(\lambda)$ containing the unit element for 
arbitrary rank 2 cases. 
We also present the explicit form of $B_0(\lambda)$ 
containing the highest (or lowest) weight vector 
by the polyhedral realization method.
\end{abstract}
\tableofcontents

%%%%%%%%%%%%%%% Section 1 %%%%%%%%%%%%%%%%
\renewcommand{\thesection}{\arabic{section}}
\section{Introduction}
\setcounter{equation}{0}
\renewcommand{\theequation}{\thesection.\arabic{equation}}

\q\,\,\,In 1985, Drinfel'd and Jimbo introduced independently 
the quantum algebra $\uq:={\lan e_i, f_i, q^h \ran}_{i \in I}$
$(I= \{1,2,\cd,n\})$, which is called $q$-$analogue$ of the
universal enveloping algebra for a 
symmetrizable Kac-Moody Lie algebra $\ge$.   
The nilpotent part $\uqm$$(={\lan f_i \ran}_{i \in I})$, 
$\uqp$$(={\lan e_i \ran}_{i \in I})$ and the integrable $\uq$-modules
have the ``crystal base'' and the crystal base of $\uqm$ (resp. 
$\uqp$) is denoted by $(L(\ify), B(\ify))$ (resp. $(L(-{\ify}), 
B(-{\ify}))$).

The theory of ``crystal base'' was established by Kashiwara 
\cite{K0}, \cite{K1} and
it developes the representation theory of $\uq$
in the side of combinatrics. We can deal with
the representation theory using combinatorial methods.
One of the great 
properties of crystal base is that a tensor product of crystal
base is again a crystal base of tensor product of corresponding 
modules \cite{K1}.
We define the ``crystal '' as a combinatorial notion abstracting 
the properties of
crystal base without assuming the corresponding modules.
We can define the tensor product structure on crystals 
in a similar manner to crystal bases.

The modified quantum algebra 
$\util:=\bigoplus_{\lm \in P}\uq a_{\lm}$ 
(resp. $\uq a_{\lm}:=\uq/\sum\uq(q^h-q^{\lan h_i,\lm\ran}))$ is defined
by modifying the Cartan part of $\uq$.
Since $\util$-modules always have weight space decompositions, $\util$ is
an algebra more appropriate than $\uq$ for the
research of the modules of the category ${\cal O}_{\rm int}$.
Lusztig \cite{L1} showed that $\util$ has a crystal base $(L(\util),B(\util))$
(resp. $(L(\uq a_{\lm}),B(\uq a_{\lm}))$) and Kashiwara \cite{K4}
described that the existence of the following isomorphism of
crystals:
$$
B(\uq a_{\lm})\cong B(\ify)\ot T_{\lm}\ot B(-{\ify})
$$
where $T_{\lm}$ is the crystal given in Section 2.

Polyhedral realizations of crystal bases is one of the methods 
for describing the crystal base explicitly, which was introduced
by Nakashima and Zelevinsky \cite{NZ}. 
We can describe a vector of crystal bases as a lattice point of 
certain convex polyhedron in an infinite $\mathbb Z$-lattice 
by this method. This method can be applied to
not only classical types but also affine or more general
Kac-Moody types.
In \cite{NZ}, polyhedral realization of $B(\ify)$ is given when
$\ge$ is of arbitrary rank 2 cases, of $A_n$ and of type $A^{(1)}_{n-1}$ and
in \cite{N1}, Nakashima gave the polyhedral realization of the crystal base 
$B(\lm)$ $(\lm \in P_+)$ of irreducible 
integrable highest weight module when $\ge$ is in the same cases as above. 
He and the author \cite{HN} applied this method 
to the modified quantum algebras
and had the polyhedral realization of $B(\uq a_{\lm})$, which is,
in general, not connected. So, we also described some specific
connected component of $B(\lm)$ $(\lm \in P)$ containing 
$u_{\ify} \ot t_{\lm} \ot u_{-{\ify}}$ (we denote this connected
component by $B_0(\lm)$) and the explicit form of 
the highest weight vector in $B_0(\lm)$ under certain assumption
on the weight $\lm$ for $\ge = A_n$ and $A^{(1)}_1$.
\vskip5pt

In this paper, we give the polyhedral realization of the
connected component $B_0(\lm)$ containing the highest (or lowest) 
weight vector for the modified quantum algebras of arbitrary rank 2 
cases. When $\ge$ is of classical type, we know that $B_0(\lm)$
contains a highest weight vector and a lowest one \cite{K4}, \cite{L1}. 
But when $\ge$ is
of affine or hyperbolic type, $B_0(\lm)$ may not necessarily
contain either a highest weight vector or a lowest one.
Hence, we give the conditions of a weight $\lm$ 
for existence of the highest (or lowest) weight
vector in $B_0(\lm)$ for affine or hyperbolic type and 
describe the explicit form of 
the highest (or lowest) weight vector in $B_0(\lm)$.
Furthermore, using the explicit form of highest (or lowest)
weight vector, we present the polyhedral realization of $B_0(\lm)$.
\vskip5pt
This paper is organized as follows:
in Section 2, we recall the definition of the modified quantum
algebra and its crystal base.
In Section 3, we review the method of polyhedral realization and give
the polyhedral realization of $B(\pm{\ify})$ and $B(\uq a_{\lm})$.
In Section 4, for the modified quantum algebra of affine or hyperbolic
type of rank 2 we describe the explicit form of the highest (or lowest)
weight vector in the connected component $B_0(\lm)$
and give the condition of a weight $\lm$ for existence 
of the highest (or lowest)
weight vector in $B_0(\lm)$. In Section 5, 
we present the polyhedral realization of $B_0(\lm)$ containing 
the highest (or lowest) weight vector.
In the appendix, we show the results of classical types of
$A_2$, $B_2$ and $G_2$.

%%%%%%%%%% section 2  %%%%%%%%%%%%%%%
\section{Modified quantum algebra and its crystal base}
\setcounter{equation}{0}
\renewcommand{\thesection}{\arabic{section}}
\renewcommand{\theequation}{\thesection.\arabic{equation}}

%%%%%%%%%%% 2.1 %%%%%%%%%%%%%
%\subsection{Modified Quantum Algebra}

In this section, we review the definition of the
modified quantum algebra $\util$ and its crystal base. 
First of all, we define the quantum algebra $\uq$.
We fix a finite index set $I$ and
let $A=(a_{ij})_{i,j\in I}$ be 
a generalized symmetrizable Cartan matrix,
$(\ttt,\{\al_i\}_{i\in I},\{h_i\}_{i\in I})$ be the associated 
Cartan data and $\ge$ be the associated Kac-Moody Lie algebra
where $\al_i$ (resp. $h_i$) is called an simple root
(resp. simple coroot).
Let $P$ be a weight lattice with a $\QQ$-valued symmetric 
bilinear form $(\; , \;)$, $P^*$ be a dual lattice 
including $\{h_i\}_{\in I}$ and 
$Q:=\bigoplus_{i\in I}\QQ(q)\al_i$ be a root lattice.
We define the quantum algebra $\uq$ to be the associative algebra with 
$1$ over $\QQ(q)$ generated by $e_i,f_i, q^h$ ($i\in I,h\in P^*$) 
with the following relations:
\begin{eqnarray*}
&q^0 = 1,\;\; q^h q^{h'}= q^{h+h'}&  \text{ for } h,h' \in P^*, \\
&q^h e_i q^{-h} = q^{\lan h, \alpha_i \ran} e_i,\;\;
q^h f_i q^{-h} = q^{-{\lan h, \alpha_i \ran}} f_i
&\text{ for } i \in I,  h \in P^*, \\
&[e_i, f_j] = \delta_{i,j}\dfrac{t_i - t_i^{-1}}
{q_i - q_i^{-1}}
&\text{ for } i,j \in I
\end{eqnarray*}
where $q_i= q^{(\alpha_i, \alpha_i)/2}$ and 
$t_i = q^{(\alpha_i, \alpha_i)h_i/2}$,
\begin{eqnarray*}
&\displaystyle{\sum^{1-a_{ij}}_{k=0}}(-1)^{k}e_i^{(k)}e_j e_i^{(1-a_{ij}-k)}=
\displaystyle{\sum^{1-a_{ij}}_{k=0}}(-1)^{k}f_i^{(k)}f_j f_i^{(1-a_{ij}-k)}=0&
\text{ for } i,j \in I, i \ne j
\end{eqnarray*}
where for $n \in \ZZ_{\geq 0}$ 
 $$[n]_i = \frac{q_i^n - q_i^{-n}}{q_i - q_i^{-1}}\,,\q
[n]_i! = \displaystyle{\prod_{k=1}^n}[k]_i $$
and we define $e_i^{(n)}= {e_i^n}/{[n]_i!}$ and $f_i^{(n)}= {f_i^n}/{[n]_i!}$.

We shall define the modified quantum algebra $\util$.
We define the left $\uq$-module $\uq a_{\lm}$ by the relation: 
$q^h a_{\lm} = q^{\lan h,\lm\ran}a_{\lm}$.
Then $\util = \oplus_{\lm \in P}\uq a_{\lm}$ has 
an algebra structure by 
\begin{enumerate}
\item
$a_{\lm}P = Pa_{\lm-\xi}$ \,\, for \,\,$\xi \in Q$
and \,\,$P \in \uq_{\xi}$\\
$(\,\uq_{\xi}
:=\{P \in \uq;\,q^hPq^{-h}=q^{\lan h,\xi\ran}P\,
{\hbox{ for any }} h \in P^*\}\,)$
\item
$a_{\lm}a_{\mu}=\del_{\lm,\mu}a_{\lm}$
\end{enumerate}
and we call this algebra {\it modified quantum algebra}.

Let ${\cal O}_{\rm int}$
be the category whose object is $\uq$-module
satisfying that it has a weight space decomposition and 
all $e_i , f_i$ $(i \in I)$ are locally nilpotent.
It is well-known that the category ${\cal O}_{\rm int}$ is 
a semisimple category and all simple ojbects are parametrized 
by dominant integral weights $P_+$. 
Let $M$ be a $\uq$-module with the weight space decomposition 
$M=\oplus_{\lm\in P}M_\lm$. Then $a_\lm$ is a projection 
$a_\lm:M\lar M_\lm$. Therefore, we think that $\util$ is
an algebra more appropriate than $\uq$ for the
research of the modules of the category ${\cal O}_{\rm int}$.
\vskip5pt

Now, we shall review the theory of crystal base.
We follow the notations and terminologies to \cite{N1}, \cite{NZ}.

Let $M$ be a $\uq$-module  in $\cal O_{\rm int}$. 
For any $u\in M_\lm$ $(\lm\in P)$,
we have the unique expression:
\[
 u=\sum_{n\geq 0}f_i^{(n)}u_n,
\]
where $u_n\in {\rm Ker}e_i\cap M_{\lm+n\al_i}$. By using this,
we define the Kashiwara operator $\eit,\fit\in {\rm End}(M)$,
\[
 \eit u:=\sum_{n\geq 1}f_i^{(n-1)}u_n,\,\,\qq
 \fit u:=\sum_{n\geq 0}f_i^{(n+1)}u_n.
\]
here note that we can define the Kashiwara operators 
$\eit,\fit\in {\rm End}(\uqpm)$ by the similar manner \cite{K1}.
Let $A\subset \QQ(q)$ be the subring of rational functions which
are regular
 at $q=0$. Let $M$ be the $\uq$-module in $\cal O_{\rm int}$.

\begin{def1}
A pair $(L,B)$ is a crystal base of $M$ $($resp. $\uqpm$$)$, 
if it satisfies following conditions:
\begin{enumerate}
\item
 $L$ is free $A$-submodule of $M$ $($resp. $\uqpm$$)$ and 
$M\cong \QQ(q)\ot_A L$ $($resp. $\uqpm\cong \QQ(q)\ot_A L$$)$.

\item
 $B$ is a basis of the $\QQ$-vector space $L/qL$.

\item
 $L=\oplus_{\lm \in P}L_{\lm}$, $B=\sqcup_{\lm \in P}B_{\lm}$ where
$L_{\lm} := L \cap M_{\lm}$, $B_{\lm} := B \cap L_{\lm}/qL$
$($resp. there is no corresponding statement$)$.

\item
 $\eit L\subset L$ and
$\fit L\subset L$.

\item
$\eit B\subset B\sqcup\{0\}$ and
$\fit B\subset B\sqcup\{0\}$
$($resp. $\fit B\subset B$$)$ $($$\eit$ and $\fit$ acts on
$L/qL$ by (iv)$)$.

\item
 For $u,v\in B$,
$\fit u=v$ if and only if $\eit v=u$.
\end{enumerate}
\end{def1}
\noindent
The unit of subalgebra $\uqmp$ is denoted by
 $u_{\pm\ify}$. We set
\beqnn
L(\pm\ify) & :=& \sum_{i_j\in I,l\geq 0}
A\til f_{i_l}\cd \til f_{i_1}u_{\pm\ify},
\\
B(\pm\ify) & := &
\{\til f_{i_l}\cd \til f_{i_1}u_{\pm\ify}\,\,
{\rm mod}\,\,qL(\ify)\,|\,i_j\in I,l\geq 0\}
\setminus \{0\},
\eeqnn
then we have
\begin{thm1}[\cite{K1}]
A pair $(L(\pm\ify),B(\pm\ify))$ 
is a crystal base of $\uqmp$.
\end{thm1}

Now we introduce the notion {\it crystal}, which is obtained by
abstracting the combinatorial properties of crystal bases.
\begin{def1}
A {\it crystal} $B$ is a set endowed with the following maps:
\begin{eqnarray*}
&& wt:B\lar P,\\
&&\vep_i:B\lar\ZZ\sqcup\{-\infty\},\q
  \vp_i:B\lar\ZZ\sqcup\{-\infty\} \q{\hbox{for}}\q i\in I,\\
&&\eit:B\sqcup\{0\}\lar B\sqcup\{0\},
\q\fit:B\sqcup\{0\}\lar B\sqcup\{0\}\q{\hbox{for}}\q i\in I,\\
&&\eit(0)=\fit(0)=0,
\end{eqnarray*}
and those maps satisfy the following axioms: for
 all $b,b_1,b_2 \in B$, we have
\begin{eqnarray*}
&&\vp_i(b)=\vep_i(b)+\lan h_i,wt(b)\ran,\\
&&wt(\eit b)=wt(b)+\al_i{\hbox{ if  }}\eit b\in B,\\
&&wt(\fit b)=wt(b)-\al_i{\hbox{ if  }}\fit b\in B,\\
&&\eit b_2=b_1 \Longleftrightarrow \fit b_1=b_2\,\,(\,b_1,b_2 \in B),\\
&&\vep_i(b)=-\ify
   \Longrightarrow \eit b=\fit b=0.
\end{eqnarray*}
\end{def1}
\noindent
Indeed, if $(L,B)$ is a crystal base, then $B$ is a crystal.

\begin{def1}
\begin{enumerate}
\item
Let $B_1$ and $B_2$ be crystals. A {\it strict morphism} of crystals
$\psi:B_1\lar B_2$ is a map
$\psi:B_1\sqcup\{0\} \lar B_2\sqcup\{0\}$
satisfying the following:
(1) $\psi(0)=0$. 

(2) If $b\in B_1$ and $\psi (b)\in B_2$, then
\[
\hspace{-30pt}wt(\psi(b))=wt(b),\q \vep_i(\psi(b))=\vep_i(b),\q
  \vp_i(\psi(b))=\vp_i(b).
\]
and the map $\psi$ commutes with all $\eit$ and $\fit$.
\item
An injective strict morphism is called an 
{\it embedding of crystals}. We
call $B_1$ is a subcrystal of $B_2$, if $B_1$ is a subset of $B_2$ and
becomes a crystal itself by restricting the data on it from $B_2$.
\end{enumerate}
\end{def1}
\noindent
The following examples of crystals will play an important role
in the subsequent sections.
\begin{ex1}
\label{ex-tlm}
Let  $T_{\lm}:=\{t_{\lm}\}$ $(\lm\in P)$
be the crystal consisting of one element $t_\lm$ defined by
$wt(t_{\lm})=\lm,$
$\vep_i(t_{\lm})=\vp_i(t_{\lm})= -{\ify}$
, $\eit (t_{\lm})=\fit(t_{\lm})=0$.
\end{ex1}

\begin{ex1}
\label{ex-bi}
For $i\in I$, the crystal $B_i:=\{(x)_i\,: \, x \in \ZZ\}$ 
is defined by:
\beqnn
&& wt((x)_i)=x \al_i,\qq \vep_i((x)_i)=-x, \qq \vp_i((x)_i)=x,\\
&& \vep_j((x)_i)=-\ify,\qq \vp_j((x)_i)=-\ify \q {\rm for }\q j\ne i,\\
&& \til e_j (x)_i=\del_{i,j}(x+1)_i,\qq
\til f_j(x)_i=\del_{i,j}(x-1)_i.
\eeqnn
\end{ex1}
\noindent
Note that as a set $B_i$ is identified with the set of integers
$\mathbb Z$.

%%%%%%%%%%% 2.2 %%%%%%%%%%%%%
%\subsection{Crystal base of $\uq a_{\lm}$}

\vskip5pt
Here we see the properties of the crystal base for the modified quantum
algebra $\util$.
Lusztig  \cite{L1} showed that 
$\util$ has a crystal base and Kashiwara  \cite{K4} described
the existence of the following isomorphism of crystals:

\begin{thm1}[\cite{K4}]
\begin{eqnarray*}
B(\uq a_{\lm}) &\cong& B(\ify)\ot T_{\lm}\ot B(-{\ify}),\\
B(\util) &\cong& \bigoplus_{\lm\in P}B(\ify)\ot T_{\lm}\ot B(-{\ify}).
\end{eqnarray*}
\end{thm1}

%%%%%%%%%% section 3  %%%%%%%%%%%%%%%
\section{Polyhedral realization of 
$B(\pm{\ify})$ and $B(\uq a_{\lm})$}
\setcounter{equation}{0}
\renewcommand{\thesection}{\arabic{section}}
\renewcommand{\theequation}{\thesection.\arabic{equation}}

%%%%%%%%%% 3.0 %%%%%%%%%%%%%
In this section, 
we review the polyhedral realization of the crystal $B(\pm\ify)$ 
(see \cite{NZ})
and $B(\uq a_{\lm})$ (see \cite{HN}). At first, we shall 
recall the polyhedral realization of $B(\pm\ify)$.
We consider the following additive groups:
\begin{eqnarray*}
\ZZ^{+\ify}&:=&\{(\cd,x_k,\cd,x_2,x_1)\,|\, x_k\in\ZZ
\,\,{\rm and}\,\,x_k=0\,\,{\rm for}\,\,k\gg 0\},\\
\ZZ^{-{\ify}}&:=&\{(x_{-1},x_{-2}, \cd , x_{-k}, \cd)\,|\, x_{-k} \in \ZZ 
\,\,{\rm and}\,\, x_{-k}=0 \,\,{\rm for} \,\,k\gg 0 \}.
\end{eqnarray*}
We will denote by
 $\ZZ^{+{\ify}}_{\geq 0} \subset \ZZ^{+{\ify}}$
(resp. $\ZZ^{-{\ify}}_{\leq 0} \subset \ZZ^{-{\ify}}$) 
the semigroup of nonnegative (resp. nonpositive) sequences. 
Take an infinite sequence of indices
$\io^+=(\cd,i_k,\cd,i_2,i_1)$ 
(resp. $\io^-=(i_{-1},i_{-2},\cd , i_{-k},\cd)$) from $I$ such that
\begin{equation}
{\hbox{
$i_k\ne i_{k+1}$ for any $k$, and $\sharp\{k>0 \,\,
{\rm (resp.\,\, k<0)}\,\,: i_k=i\}=\ify$ for any $i\in I$.}}
\label{seq-con}
\end{equation}
We can have a crystal structure on $\ZZ^{+{\ify}}$ 
(resp. $\ZZ^{-{\ify}}$) associate to 
 $\io^+$ (resp. $\io^-$) (see \cite{NZ}) and denote it by
 $\ZZ^{+{\ify}}_{\io^+}$ 
(resp. $\ZZ^{-{\ify}}_{\io^-}$). 
Let $B_i$ be the crystal given in Example \ref{ex-bi}. We 
obtain the following embeddings (\cite{K2}):
\begin{eqnarray*}
\Psi_{i}^+ :& B(\infty)& \hookrightarrow \; B(\infty)\ot B_i
\q
(u_{\ify}\mapsto u_{\ify}\ot (0)_i),\\
\Psi_{i}^- :& B(-{\infty})& \hookrightarrow \; B_i \ot B(-{\infty})
\q
(u_{-{\ify}}\mapsto (0)_i \ot u_{-{\ify}}).
\end{eqnarray*}
Iterating $\Psi_{i}^+$ (resp. $\Psi_{i}^-$)
according to $\io^+$ (resp. $\io^-$), we get 
the {\it Kashiwara embedding} (\cite{K2}):
\begin{eqnarray}
\Psi_{\io^+}:&B(\ify)& \hookrightarrow \; \ZZ^{+{\ify}}_{\geq 0}
\subset \ZZ^{+{\ify}}_{\io^+}
\q
(u_{{\ify}}\mapsto (\cd,0,\cd,0,0,0)),
\label{kas+}\\
\Psi_{\io^-}:&B(-{\ify})&\hookrightarrow \; \ZZ^{-{\ify}}_{\leq 0}
\subset \ZZ^{-{\ify}}_{\io^-}
\q
(u_{-{\ify}}\mapsto (0,0,0,\cd,0,\cd)).
\label{kas-}
\end{eqnarray}
We consider the following infinite dimensional vector spaces
and their dual spaces:
\begin{eqnarray*}
\QQ^{+{\ify}}&:=&\{\vec{x}=
(\cd,x_k,\cd,x_2,x_1): x_k \in \QQ\,\,{\rm and }\,\,
x_k = 0\,\,{\rm for}\,\, k \gg 0\},\\
\QQ^{-{\ify}}&:=&\{\vec{x}=
(x_{-1},x_{-2},\cd,x_{-k},\cd): x_{-k} \in \QQ\,\,{\rm and }\,\,
x_{-k} = 0\,\,{\rm for}\,\, k \gg 0\},\\
 (\QQ^{\pm{\ify}})^*&:=&{\rm Hom}(\QQ^{\pm{\ify}},\QQ).
\end{eqnarray*}
We will write a linear form $\vp \in (\QQ^{+{\ify}})^*$ as
$\vp(\vec{x})=\sum_{k \geq 1} \vp_k x_k$ ($\vp_j\in \QQ$). 
Similarly, we write 
$\vp \in (\QQ^{-{\ify}})^*$ as
$\vp(\vec{x})=\sum_{k \leq -1} \vp_k x_k$ ($\vp_j\in \QQ$).

For the sequence 
$\io^+=(i_k)_{k \geq 1}$ (resp. $\io^-=(i_k)_{k \leq -1}$) and 
$k \geq 1$ (resp. $k \leq -1$), we set 
\[
\kp:={\rm min}\{l:l>k>0\,\, ({\rm resp.}\,0>l>k)\,\,
{\rm  and }\,\,i_k=i_l\},
\]
if it exists, and
$$
\km:={\rm max}\{l:0<l<k\,\,({\rm resp.}\,l<k<0)\,\,
{\rm  and }\,\,i_k=i_l\},
$$
if it exists, otherwise $\kp=\km=0$.

We define a linear form $\beta_k$ $(k\geq0)$ on $\QQ^{+{\ify}}$ by
\begin{equation}
\beta_{k}(\vec x) := \begin{cases}
x_k+\sum_{k<j<\kp}\lan h_{i_k},\al_{i_j}\ran x_j+x_{\kp}
&
(k \geq 1),\\
0& (k = 0).
\end{cases}
\end{equation}
We also define a linear form $\beta_k$ $(k\leq 0)$ 
on $\QQ^{-{\ify}}$ by 
\begin{equation}
\beta_k(\vec y):= \begin{cases}
y_{\km}+\sum_{\km<j<k}\lan h_{i_k},\al_{i_j}\ran y_j+y_k
& 
(k \leq -1),\\
0 &(k=0).
\end{cases}
\end{equation}
By using these linear forms, let us
 define a piecewise-linear operator 
$S_k=S_{k,\io}$ on $(\QQ^{\pm{\ify}})^*$
as follows:
\begin{equation}
S_k(\vp):=\begin{cases}
\vp-\vp_k\beta_k & \text{if}\;\, \vp_k>0,\\
 \vp-\vp_k\beta_{\km} & \text{if}\;\, \vp_k\leq 0,\\
\end{cases}
\label{Sk}
\end{equation}
for $\vp(\vec x)=\sum \vp_k x_k\in (\QQ^{\pm{\ify}})^*$.
Here we set
\beqnn
\Xi_{\io^\pm} &:= &
\{S_{\pm j_l}\cd S_{\pm j_2}S_{\pm j_1}(\pm x_{\pm{j_0}})\,|\,
l\geq0,j_0,j_1,\cd,j_l\geq1\},\\
\Sigma_{\io^\pm} & := &
\{\vec x\in \ZZ^{\pm{\ify}}\subset \QQ^{\pm{\ify}}
\,|\,\vp(\vec x)\geq0\,\,{\rm for}\,\,
{\rm any}\,\,\vp\in \Xi_{\io^\pm}\}.
\eeqnn
We impose on $\io^+$ and $\io^-$ the following assumptions (P),(N):
\begin{eqnarray}
&{\hbox{(P)\;   for $\io^+$, if $\km=0$ then $\vp_k\geq0$ for any 
$\vp(\vec x)=\sum_k\vp_kx_k\in \Xi_{\io^+}$}},\label{p}\\
&{\hbox{(N)\;   for $\io^-$, if $\kp=0$ then $\vp_k\leq0$ for any 
$\vp(\vec x)=\sum_k\vp_kx_k\in \Xi_{\io^-}$}}.\label{n}
\end{eqnarray}

\begin{thm1}[\cite{NZ}]
Let $\io^{\pm}$ be the indices of sequences which are satisfied 
$(\ref{seq-con})$ and the assumptions $(P)$, $(N)$. 
Suppose $\Psi_{\io^+}:B(\ify)\hookrightarrow \ZZ^{\ify}_{\io^+}$ and
$\Psi_{\io^-}:B(-{\ify})\hookrightarrow \ZZ^{-{\ify}}_{\io^-}$ 
are the 
Kashiwara embeddings. Then, we have
${\rm Im}(\Psi_{\io^+})(\cong B(\ify))=\Sigma_{\io^+}$,\,
${\rm Im}(\Psi_{\io^-})(\cong B(-{\ify}))=\Sigma_{\io^-}$.
\end{thm1}
\noindent
We call $\Sigma_{\io^+}$ (resp. $\Sigma_{\io^-}$) the polyhedral
realization of $B(\ify)$ (resp. $B(-{\ify})$).

\vskip5pt

%%%%%%%%%%% 3.1 %%%%%%%%%%%%%
%\subsection{Crystal structure of $\ZZ^{\ify}[\lm]$}

Here, we will recall the polyhedral realization of 
$B(\uq a_{\lm})$.
As a set, we denote  
$\ZZ^{+{\ify}} \ot T_{\lm} \ot \ZZ^{-{\ify}}$
(resp. $\QQ^{+{\ify}} \ot T_{\lm} \ot \QQ^{-{\ify}}$)
by $\ZZ^{\ify}[\lm]$ (resp. $\QQ^{\ify}[\lm]$).
We shall construct the crystal structure of 
$\ZZ^{\ify}[\lm]$.
We fix the indices sequence $\io:=(\io^+,t_{\lm},\io^-)=
(\cd,i_2,i_1,t_{\lm},i_{-1},i_{-2},\cd)$ and a weight $\lm\in P$.
Since $\ZZ^{\ify}[\lm]$ is the subset of $\QQ^{\ify}[\lm]$, we
can denote $\vec x \in \ZZ^{\ify}[\lm]$ as 
$\vec x=(\cd,x_2,x_1,t_{\lm},x_{-1},x_{-2},\cd)$.

%For $\ZZ^{+{\ify}}_{{\io}^+}$, $\ZZ^{-{\ify}}_{{\io}^-}$ in section 2.2, 
%we denote 
%$\ZZ^{\ify}_{\io}[\lm]$ to be 
%$\ZZ^{\ify}_{{\io}^+}\ot T_{\lm}\ot \ZZ^{-{\ify}}_{{\io}^-}.
%$
%We will construct the crystal structure of 
%$\ZZ^{\ify}[\lm]$.
%We fix the indices sequence $\io:=(\io^+,t_{\lm},\io^-)=
%(\cd,i_2,i_1,t_{\lm},i_{-1},i_{-2},\cd)$ and weight $\lm\in P$.
%
%We identify $\ZZ^{+{\ify}} \ot T_{\lm} \ot \ZZ^{-{\ify}}$ with
%$\ZZ^{\ify}[\lm]$. Therefore, $\ZZ^{\ify}[\lm]$ is the subset of
%$\QQ^{\ify}$. Thus, we can deonte $\vec x \in \ZZ^{\ify}[\lm]$ as 
%$\vec x=(\cd,x_2,x_1,t_{\lm},x_{-1},x_{-2},\cd)$.
For $\vec x \in \QQ^{\ify}[\lm]$, we define a linear function 
$\sigma_k(\vec x)\,\,(k \in \ZZ)$ by:
\beqn
\sigma_k(\vec x):= 
  \begin{cases}
    x_k+\sum_{j>k}\lan h_{i_k},\al_{i_j}\ran x_j
   &(k\geq1),\\
    -\lan h_{i_k},\lm\ran+x_k+\sum_{j>k}\lan h_{i_k},\al_{i_j}\ran x_j
   &(k\leq -1),\\
    -\ify
   &(k=0).
\end{cases}
\label{sigma}
\eeqn
Since $x_j=0$ for $j\gg0$,
 $\sigma_k$ is
well-defined.
Let $\sigma^{(i)} (\vec x)
 := {\rm max}_{k: i_k = i}\sigma_k (\vec x)$ and
\begin{equation}
M^{(i)} = M^{(i)} (\vec x) :=
\{k: i_k = i, \sigma_k (\vec x) = \sigma^{(i)}(\vec x)\}.
\label{m(i)}
\end{equation}
Note that 
$\sigma^{(i)} (\vec x)\geq 0$, and that 
$M^{(i)} = M^{(i)} (\vec x)$ is finite set if and only if 
$\sigma^{(i)} (\vec x) > 0$. 
Now, we define the map 
$\eit: \ZZ^{\ify}[\lm] \sqcup\{0\}\lar \ZZ^{\ify}[\lm] \sqcup\{0\}$
, 
$\fit: \ZZ^{\ify}[\lm] \sqcup\{0\}\lar \ZZ^{\ify}[\lm] \sqcup\{0\}$ , by
 $\eit(0)=\fit(0)=0$ and 
\begin{eqnarray}
\hspace{-10pt}(\fit(\vec x))_k &=& x_k + \delta_{k,{\rm min}\,M^{(i)}}
\,\,{\rm if }\,\,M^{(i)} \,\,{\rm exists};
\,\,{\rm otherwise}\,\,\fit(\vec x)=0,
\label{action-f}\\
\hspace{-10pt}(\eit(\vec x))_k &=& 
x_k - \delta_{k,{\rm max}\,M^{(i)}} \,\, {\rm if}\,\,
M^{(i)} \,\,{\rm exists}; \,\,
 {\rm otherwise} \,\, \eit(\vec x)=0.
\label{action-e}
\end{eqnarray}
where $\del_{i,j}$ is Kronecker's delta.
We also define the weight function and the function
$\vep_i$ and $\vp_i$ on $\ZZ^{\ify}[\lm]$ as follows:
\begin{equation}
\begin{array}{l}
wt(\vec x) :=\lm -\sum_{j=-{\ify}}^{\ify} x_j \al_{i_j}, \,\,
\vep_i (\vec x) := \sigma^{(i)} (\vec x),\\
\vp_i (\vec x) := \lan h_i, wt(\vec x) \ran + \vep_i(\vec x).
\end{array}
\label{wt-vep-vp}
\end{equation}
We denote this crystal by $\ZZ^{\ify}_{\io}[\lm]$.
We can identify $\ZZ^{\ify}_{\io}[\lm]$ with 
$\ZZ^{\ify}_{{\io}^+}\ot T_{\lm}\ot \ZZ^{-{\ify}}_{{\io}^-}$.

Since there exist the embedings of crystals:
$B(\pm\ify)\hookrightarrow \ZZ_{\io^\pm}^{\pm\ify}$, we obtain
\begin{thm1}[\cite{HN}]
\begin{eqnarray*}
\Psi_{\io}^{(\lm)}:B(\ify)\ot T_{\lm}\ot B(-{\ify})
&\hookrightarrow&
\ZZ^{+{\ify}}_{{\io}^+}\ot T_{\lm}\ot \ZZ^{-{\ify}}_{{\io}^-}(=\ZZ^{\ify}_{\io}[\lm])\\
u_{\ify}\ot t_{\lm}\ot u_{-{\ify}}\q &\mapsto& \q(\cd,0,0,t_{\lm},0,0,\cd)
\end{eqnarray*}
is the unique stirict embedding which is associated with 
$\io:=(\cd,i_2,i_1,t_{\lm},i_{-1},i_{-2},\cd)$.
\label{Psi-lm}
 \end{thm1}

%%%%%%%%%%%% 3.2 %%%%%%%%%%%%%%%%
Now, we will give the polyhedral realization of $B(\uq a_{\lm})$.
We fix an indices sequence $\io$.
We define a linear function $\bar\beta_k(\vec x)$ as follows:
\beqn
\bar\beta_k (\vec x) & = & \sigma_k (\vec x) - \sigma_{\kp} (\vec x)
\label{beta}
\eeqn
where $\sigma_k$ is defined by (\ref{sigma}).
Since $\lan h_{i},\al_{i}\ran = 2$ for any $i \in I$, we have 
\beqnn
\bar\beta_k(\vec{ x}) & = &
  \begin{cases}
    x_k+\sum_{k<j<\kp}
    \lan h_{i_k},\al_{i_j}\ran x_j+x_{\kp}
    &(k\geq1\,\,\text{or}\,\,\kp \leq -1),\\
    -\lan h_{i_k},\lm\ran + 
    x_k+\sum_{k<j<\kp}
    \lan h_{i_k},\al_{i_j}\ran x_j+x_{\kp}
    &(k\leq-1\,\,\text{and}\,\,\kp > 0).
\end{cases}
\eeqnn
Using this notation, we define an operator $\bar{S_k}=\bar{S}_{k,\io}$
for a linear function $\vp(\vec x)=c+\sum_{-{\ify}}^{\ify}\vp_kx_k$
$(c,\vp_k\in\QQ)$ as follows:
\begin{equation}
 \bar{S_k}\,(\vp) :=\begin{cases}
 \vp - \vp_k \bar\beta_k & \text{if}\;\, \vp_k > 0,\\
 \vp - \vp_k \bar\beta_{k^{(-)}} & \text{if}\;\, \vp_k \leq 0.
\end{cases}
\label{S_k}
\end{equation}
An easy check shows $(\bar{S_k})^2=\bar{S_k}$.
For a sequence $\io$ and an integral weight $\lm$, we denote by
$\Xi_{\io}[\lm]$ the subset of linear forms that are obtained from 
the coordinate forms $x_j$, $x_{-j}$ ($j\geq1$) by applying 
transformations $\bar{S}_k$. In other words, we set
 \begin{align}
  \Xi_{\io}[\lm]:&=\{\bar{S}_{j_l}\cd \bar{S}_{j_1}(x_{j_0})\,
  :\,l\geq0,\,j_0,\cd,j_l\geq 1\} \nn
  \\
  &\cup\{\bar{S}_{{-j}_k}\cd \bar{S}_{{-j}_1}(-x_{{-j}_0})\,
  :\,k\geq0,\,j_0,\cd,j_k \geq 1\}, \label{Xi}\\
 \Sigma_{\io}[\lm]
 :&=\{\vec x\in \ZZ^{\ify}_{\io}[\lm](\subset \QQ^{\ify}[\lm])\,:\,
\vp(\vec x)\geq 0\,\,{\rm for \,\,any }\,\,\vp\in \Xi_{\io}[\lm]\}.
\label{Sigma}
\end{align}

%\begin{equation}
% \begin{array}{ll}
%  \Xi_{\io}[\lm]&:=\{\bar{S}_{j_l}\cd \bar{S}_{j_1}(x_{j_0})\,
%  :\,l\geq0,\,j_0,\cd,j_l\geq 1\}
%  \\
%  &\cup\{\bar{S}_{{-j}_k}\cd \bar{S}_{{-j}_1}(-x_{{-j}_0})\,
%  :\,k\geq0,\,j_0,\cd,j_k\geq1\}.
% \end{array}
%\label{Xi}
%\end{equation}
%
%Now we set 
%\begin{equation}
% \Sigma_{\io}[\lm]
% :=\{\vec x\in \ZZ^{\ify}_{\io}[\lm](\subset \QQ^{\ify})\,:\,
%\vp(\vec x)\geq 0\,\,{\rm for \,\,any }\,\,\vp\in \Xi_{\io}[\lm]\}.
% \label{Sigma}
%\end{equation}

\begin{thm1}[\cite{HN}]
 \label{main}
Suppose that $\io^{\pm}$ is satisfied $(\ref {seq-con})$
and the assumptions $(P)$, $(N)$. 
Let $\Psi^{(\lm)}_{\io}:
B(\ify)\ot T_{\lm}\ot B(-{\ify})\hookrightarrow 
\ZZ^{\ify}_{\io}[\lm]$
be the embedding of Theorem $\ref{Psi-lm}$. Then, 
${\rm Im}(\plm)(\cong B(\ify)\ot T_{\lm}\ot B(-{\ify}))$ 
$=$ 
$\Sigma_{\io}[\lm]$.
\end{thm1}
\noindent
We call $\Sigma_{\io}[\lm]$ the polyhedral realization of
$B(\uq a_{\lm})$.

The following lemma will be used in Section 5.

\begin{lem1}[\cite{HN}]
\label{mlem}
Let $\Xi_{\io}$ be a set of linear functions on $\QQ^{\ify}[\lm]$.
Suppose $\Xi_{\io}$ 
is closed by actions of $\bar{S}_k$, then 
the set 
\[
\Sigma_{\io}:=\{\vec x\in \ZZ^{\ify}_{\io}[\lm]\,:\,
\vp(\vec x)\geq0{\text{ for any }}\vp\in \Xi_\io\}
\] 
holds a crystal structure.
 \end{lem1}

%%%%%%%%% Section 4  %%%%%%%%%%%%%
\section{Highest and lowest weight vector of $B_0(\lm)$ for affine and 
hyperbolic types of rank 2}
\setcounter{equation}{0}
\renewcommand{\theequation}{\thesection.\arabic{equation}}
%%%%%%%%%%%%%%%%%%%%

In this section, we consider the case that $\ge$ is of affine or 
hyperbolic type of rank 2. 
We define that $B_0(\lm)$ is a connected component of
$B(\uq a_{\lm})$ $(\cong B(\ify)\ot T_{\lm}\ot B(-{\ify}))$ containing
$u_{\ify} \ot t_{\lm} \ot u_{-{\ify}}$.
We fix:

\begin{itemize}
 \item$ $Cartan Matrix : $ \begin{pmatrix}
         2&-c_1\\
         -c_2&2
        \end{pmatrix}$ \; ($c_1c_2\geq4$).
 \item$\lm=\lm_1\Lm_1 + \lm_2\Lm_2$ $ \in P$ \; 
  ($\lm_1 \in \ZZ_{>0},\; \lm_2 \in \ZZ_{<0},\; \Lm_i$: fundamental weight).
 \item$\io = (\cd,1,2,1,t_{\lm},2,1,2,\cd)$ ($I=\{1,2\}$).
\end{itemize}
Note that if $\lm$ is a dominant (resp. antidominant) weight, we know
the highest (resp. lowest) weight vector of $B_0(\lm)$ is 
$(\cd,0,0,t_{\lm},0,0,\cd)$ in $\ZZ^{\ify}_{\io}[\lm]$. 
So, we consider the case that
a weight $\lm$ is not a dominant or antidominant weight.
Here, we will describe the explicit form of 
the highest (or lowest) weight vectors. 
For the purpose, we need some preparations.  
We define $Chebyshev$ $polynomials$ 
and give some facts.

\begin{def1}
   $Chebyshev$ $polynomials (:=P_k(X))$ is given by:
  \begin{equation*}
   P_k(\alpha + {\alpha}^{-1}) = \frac{{\alpha}^{k+1}-{\alpha}^{-k-1}}
                                      {\alpha - {\alpha}^{-1}}.
  \end{equation*}
 \end{def1}

\begin{fac1}
  \begin{enumerate}
   \item
     The generating function for $P_k(X)$ is given by:
$$\sum_{k\geq0}P_k(X)z^{k} = (1 - Xz + z^{2})^{-1}.$$
   \item
     $P_k(X)$ satisfies the following identities:
   \begin{align}
&(X+2)P_k(X)^{2}-(P_{k+1}(X)+P_{k}(X))(P_{k}(X)+P_{k-1}(X))=1,
 \label{id1}\\
&(P_{k}(X)+P_{k-1}(X))^{2}-(X+2)P_{k}(X)P_{k-1}(X)=1.
 \label{id2}
   \end{align}
   \item
     $P_k(X)$ satisfies the following recusion:
 $$P_k(X)=X P_{k-1}(X)-P_{k-2}(X).$$
  \end{enumerate}
 \label{fact}
 \end{fac1}

\begin{def1}
For $X=c_1c_2 -2$, we define the integers $a_{l}$, ${a'_{l}}$ $(l\geq 0)$ 
as follows $($see \text{$\cite{NZ}$}$)$:
 \begin{eqnarray*}
  &a_0=a'_0=0,\; a_1=a'_1=1,\\
  &a_{2k}=c_{1}P_{k-1}(X),\,\, {a'_{2k}}=c_{2}P_{k-1}(X)\q(k\geq1), \\
  &a_{2k+1}={a'_{2k+1}}=P_{k}(X)+P_{k-1}(X)\q(k\geq1).
 \end{eqnarray*}
\label{ak}
\end{def1}

\begin{ex1} Several first terms of $a_l$ are given by:
\begin{eqnarray*}
&a_1=1,\,\,a_2=c_1,\,\,a_3=c_1c_2-1,\,\,a_4=c_1(c_1c_2-2),
\,\,a_5=(c_1c_2-1)(c_1c_2-2)-1,&\\
&a'_1=1,\,\,a'_2=c_2,\,\,a'_3=c_1c_2-1,\,\,a'_4=c_2(c_1c_2-2),
\,\,a'_5=(c_1c_2-1)(c_1c_2-2)-1,&\\
&a_6=c_1(c_1c_2-1)(c_1c_2-3),\,\,a_7=c_1c_2(c_1c_2-2)(c_1c_2-3),&\\
&a'_6=c_2(c_1c_2-1)(c_1c_2-3),\,\,a'_7=c_1c_2(c_1c_2-2)(c_1c_2-3).&
\end{eqnarray*}
\label{exa}
\end{ex1}
\noindent
Note that if $\ge$ is of affine or hyperbolic type, $a_{l}$, $a'_{l}$
are positive.

\begin{fac1}
 $a_l$ satisfies following recusions $($see \text{$\cite{NZ}$}$)$:
 \begin{enumerate}
  \item $a_{2k+2}=c_1a_{2k+1}-a_{2k}$. 
  \item $a_{2k+1}=c_2a_{2k}-a_{2k-1}$.
  \item $a'_{2k+2}=c_2a_{2k+1}-a'_{2k}$.
  \item $a_{2k+1}=c_1a'_{2k}-a_{2k-1}$.
 \end{enumerate}
\label{af1}
\end{fac1}
{\sl Proof.}\,\,
(i) We apply Definition $\ref{ak}$ to Fact \ref{fact}(iii).  
(ii) is an easy consequence of Fact \ref{fact}(i). 
(iii), (iv) is given by Definition $\ref{ak}$ and (i),(ii).
\qed

 \begin{lem1} For any $l \geq 1$,
  \begin{equation}
   \frac{a'_{l+1}}{a'_{l}} > \frac{a_{l+2}}{a_{l+1}}\,\,,\,\,\,\,\,\,\,
   \frac{a_{l+1}}{a_{l}} > \frac{a'_{l+2}}{a'_{l+1}}.
  \end{equation}
\end{lem1}
{\sl Proof.}\,\,
We have
\begin{eqnarray*}
 P_k(X)=\frac{a_{2k+2}}{c_1}=\frac{a'_{2k+2}}{c_2},\q X+2=c_1c_2
\label{pk}
\end{eqnarray*}
by the definition. Applying Definition \ref{ak} to $(\ref{id1})$, we get
\begin{eqnarray*}
&&c_1c_2\left(\frac{a_{2k+2}}{c_1}\right)\left(\frac{a'_{2k+2}}{c_2}\right)
-(a_{2k+3}a'_{2k+1})=1,\\
&&c_1c_2\left(\frac{a_{2k+2}}{c_1}\right)\left(\frac{a'_{2k+2}}{c_2}\right)
-(a'_{2k+3}a_{2k+1})=1.
\end{eqnarray*}
Then,
\begin{align}
a_{2k+2}a'_{2k+2}-a_{2k+3}a'_{2k+1}=1, \label{a1}\\
a_{2k+2}a'_{2k+2}-a'_{2k+3}a_{2k+1}=1.\label{b1}
\end{align}
Similarly, applying Definition \ref{ak} to $(\ref{id2})$ , we get
\begin{eqnarray*}
(a_{2k+1}a'_{2k+1})-
c_1c_2\left(\frac{a_{2k+2}}{c_1}\right)\left(\frac{a'_{2k}}{c_2}\right)=1,\\
(a_{2k+1}a'_{2k+1})-
c_1c_2\left(\frac{a'_{2k+2}}{c_1}\right)\left(\frac{a_{2k}}{c_2}\right)=1.
\end{eqnarray*}
Then,
\begin{align}
a_{2k+1}a'_{2k+1}-a_{2k+2}a'_{2k}=1,\label{c1}\\
a_{2k+1}a'_{2k+1}-a'_{2k+2}a_{2k}=1.
\label{d1}
\end{align}
Using $(\ref{a1})$ and $(\ref{c1})$, we get for $l \geq 1$
\begin{eqnarray*}
a_{l+1}a'_{l+1}-a_{l+2}a'_{l}=1>0.
\end{eqnarray*}
Now, since we consider the case that $\ge$ is of the affine or
hyperbolic type,
 $a_l,\,a'_l>0$ and we have
\begin{eqnarray*}
\frac{a'_{l+1}}{a'_{l}} > \frac{a_{l+2}}{a_{l+1}}.
\end{eqnarray*}
Similarly, using $(\ref{b1})$ and $(\ref{d1})$ we have
\begin{eqnarray*}
\frac{a_{l+1}}{a_{l}} > \frac{a'_{l+2}}{a'_{l+1}}.
\end{eqnarray*}
\qed
%\hfill {\sl Q.E.D.}
 \begin{cor1}
For $k \rightarrow \ify$, following sequences are converge to real
numbers $\alpha$, $\beta$\,\,$($\,$\alpha := 
\dfrac{c_1c_2+\sqrt{c_1^2c_2^2-4c_1c_2}}{2c_2}$,
 $\beta := \dfrac{c_1c_2-\sqrt{c_1^2c_2^2-4c_1c_2}}{2c_2}$ $)$.
  \begin{enumerate}
\item \q$ \dfrac{a_{2}}{a_{1}} > \dfrac{a'_{3}}{a'_{2}}>
   \dfrac{a_{4}}{a_{3}} > \dfrac{a'_{5}}{a'_{4}}>
\cd> \dfrac{a_{2k}}{a_{2k-1}} > \dfrac{a'_{2k+1}}{a'_{2k}}>
   \cd\longrightarrow \alpha.$\\
\item \q$ \dfrac{a'_{1}}{a'_{2}} < \dfrac{a_{2}}{a_{3}}<
   \dfrac{a'_{3}}{a'_{4}} < \dfrac{a_{4}}{a_{5}}
\cd < \dfrac{a'_{2k-1}}{a'_{2k}} < \dfrac{a_{2k}}{a_{2k+1}}<
   \cd\longrightarrow \beta.$
  \end{enumerate}

%$($ where $\alpha = \frac{c_1c_2+\sqrt{c_1^2c_2^2-4c_1c_2}}{2c_2}$,
% $\beta = \frac{c_1c_2-\sqrt{c_1^2c_2^2-4c_1c_2}}{2c_2}$ $)$
\label{rem1} 
\end{cor1}

 \begin{lem1}
Suppose
$v_{-l}:=(\cd,0,0,t_{\lm},x_{-1},x_{-2},\cd ,x_{-l},0,0,\cd) \in
\ZZ_{\io}[\lm]$ for any $l \geq 0$ 
$(\, v_0:=(\cd,0,0,t_{\lm},0,0,\cd)\, )$.
Then,
  \begin{enumerate}
   \item For $l \ne 0$ and any $k \geq 1$,  
$$\sigma_{-l}(v_{-l}) \geq  \sigma_{-l-2}(v_{-l})=\sigma_{-l-2k}(v_{-1}).$$ 
   \item For any $s, t \geq1$, 
$$\sigma_{-l-1}(v_{-l})=\sigma_{-l-3}(v_{-l})=\cd=\sigma_{-l-2s+1}(v_{-l}),$$
$$\sigma_t(v_{-l})=0.$$
   \item For any $k\geq1\,(\, 2k-1<l \,)$,
if \,$\sigma_{-l+2k-1}(v_{-l})=0$\,
    and 
$\sigma_{-l-1}(v_{-l})=N>0$, then
   \begin{align*}
   &(\tilde{e}_{i_{-l-1}})^N v_{-l}\;\;\,\, =
   (\cd,0,0,t_{\lm},x_{-1},x_{-2},\cd ,x_{-l},-N,0,\cd),\\
   &(\tilde{e}_{i_{-l-1}})^{N+1} v_{-l}=0.
   \end{align*}
 \end{enumerate}
\label{lem1}
\end{lem1}
{\sl Proof.}\,\,
(i) By the definition of $\sigma$, we have
$$\sigma_{-l-2}(v_{-l})=x_{-l}+\sigma_{-l}(v_{-l})\leq \sigma_{-l}(v_{-l}).$$
(ii) It is obvious by the definition of $\sigma$.\\
\q\\
(iii) For any $n$ ($1\leq n \leq N$), we shall show that
$$(\tilde{e}_{i_{-l-1}})^n v_{-l}= 
(\cd,0,0,t_{\lm},x_{-1},x_{-2},\cd ,x_{-l},-n,0,\cd)$$
by the induction on $n$.\\

Case I) $n=1$.

\noindent
By (ii), we have for any $s ,t \geq 1$, 
$$\sigma_{-l-1}(v_{-l})=\sigma_{-l-3}(v_{-l})
=\cd=\sigma_{-l-2s+1}(v_{-l})=N,$$
$$\sigma_t(v_{-l})=0.$$
By the assumption of (iii),
for any $k\geq1\,(\,2k-1<l\,)$ we have
$$\sigma_{-l+1}(v_{-l})=\sigma_{-l+3}(v_{-l})=
\cd=\sigma_{-l+2k-1}(v_{-l})=0.$$
This shows that
$$(\tilde{e}_{i_{-l-1}})v_{-l}= 
(\cd,0,0,t_{\lm},x_{-1},x_{-2},\cd ,x_{-l},-1,0,\cd).$$

Case II) $n>1$.

\noindent
We assume that
$$(\tilde{e}_{i_{-l-1}})^{n-1} v_{-l}= 
(\cd,0,0,t_{\lm},x_{-1},x_{-2},\cd ,x_{-l},-n+1,0,\cd)=:v'_{-l}.$$
Since $x_{-1},x_{-2},\cd,x_{-l}$ does not
change, we have
$$\sigma_{-l-1}(v'_{-l})=-n+1+\sigma_{-l-1}(v_{-l}) = -n+1+N > 0$$
and by (i), we have for any $k \geq 2$
$$\sigma_{-l-1}(v'_{-l}) \geq \sigma_{-l-2k+1}(v'_{-l}).$$
For any $t \geq 1$, 
$$\sigma_t(v_{-l})=0$$
and we have for any $k\geq1\,(\,2k-1<l \,)$,
$$\sigma_{-l+1}(v'_{-l})=\sigma_{-l+3}(v'_{-l})=
\cd=\sigma_{-l+2k-1}(v'_{-l})=0$$
by the assumption of (iii).
Therefore, we get
$$(\tilde{e}_{i_{-l-1}}) v'_{-l}=(\tilde{e}_{i_{-l-1}})^{n} v_{-l}= 
(\cd,0,0,t_{\lm},x_{-1},x_{-2},\cd ,x_{-l},-n,0,\cd).$$
In particular, we consider the case n=N,
$$(\tilde{e}_{i_{-l-1}})^{N} v_{-l}= 
(\cd,0,0,t_{\lm},x_{-1},x_{-2},\cd ,x_{-l},-N,0,\cd)=:v''_{-l}.$$
In this case for any $k\geq 2$,
$$0=\sigma_{-l-1}(v''_{-l})=\sigma_{-l-2k+1}(v''_{-l})+N 
>\sigma_{-l-2k+1}(v''_{-l}).$$
We conclude
$$(\tilde{e}_{i_{-l-1}})v''_{-l}=(\tilde{e}_{i_{-l-1}})^{N+1} v_{-l}=0.$$
\qed

\noindent
Now, we define for any $j,k \geq 1$
\begin{align*}
&h_{-2k+1}\;\,:=\;a'_{2k-2}{\lm_1}\,+\,a'_{2k-1}{\lm_2},\\
&h_{-2k}\q\;\,:=\;a_{2k-1}{\lm_1}\,+\,a_{2k}{\lm_2},\\ 
&H_{-j}\q\;\;:=
(\cd,0,0,t_{\lm},\, h_{-1},\,h_{-2},\cd,h_{-j+1},\,h_{-j},0,0,\cd).
\end{align*}
 \begin{rem1}
If the condition of a weight $\lm$ is 
$\;\dfrac{a'_{2k-1}}{a'_{2k-2}}>\dfrac{\lm_1}{-\lm_2}
           \geq \dfrac{a_{2k}}{a_{2k-1}}$
(resp. $\dfrac{a_{2k}}{a_{2k-1}}>\dfrac{\lm_1}{-\lm_2}
           \geq \dfrac{a'_{2k+1}}{a'_{2k}}$ ), then
$h_{-1},\,h_{-2},\,\cd,h_{-2k+1}< 0$ and $h_{-2k},\,h_{-2k-1},\cd>0$ (resp.
 $h_{-1},\,h_{-2},\,\cd,h_{-2k}< 0$ and $h_{-2k-1},\,h_{-2k-2},\cd>0$) by
Corollary $\ref{rem1}(i)$.
\label{rem15} 
\end{rem1}
 \begin{thm1}[Highest weight vector of $B_0(\lm)$]
If the condition of a weight $\lm$ is 
\,$\dfrac{a'_{2k-1}}{a'_{2k-2}}>\dfrac{\lm_1}{-\lm_2}
           \geq \dfrac{a_{2k}}{a_{2k-1}}$ $($resp. 
$\dfrac{a_{2k}}{a_{2k-1}}>\dfrac{\lm_1}{-\lm_2}
           \geq \dfrac{a'_{2k+1}}{a'_{2k}}$$)$, then
$H_{-2k+1}$ $($resp. $H_{-2k}$$)$ is the
   highest weight vector of $B_0(\lm)$.
$($If \,$a'_{2k-2}=0$, we define that the condition of a weight 
$\lm$ is $\dfrac{\lm_1}{-\lm_2}
           \geq \dfrac{a_{2}}{a_{1}}$.$)$
\label{hwv}
 \end{thm1}
We shall prove the former case 
$\dfrac{a'_{2k-1}}{a'_{2k-2}}>\dfrac{\lm_1}{-\lm_2}
\geq \dfrac{a_{2k}}{a_{2k-1}}$ since the latter case is 
proved by the same argument of the former case.
We prepare some lemmas. 

\begin{lem1} 
We give some properties of $H_{-2k+1}$, $H_{-2k}$:
  \begin{enumerate}
  \item For any $j$ $(1 \leq j \leq 2k-1)$, $\sigma_{-j}(H_{-2k+1})=0$.
  \item For any $j$ $(1 \leq j \leq 2k)$, $\sigma_{-j}(H_{-2k})=0$.
  \item $\sigma_{-2k}(H_{-2k+1})=-a_{2k-1}{\lm_1}-a_{2k}{\lm_2}\;
(=-h_{-2k})$.
  \item
 $\sigma_{-2k-1}(H_{-2k})=-a'_{2k}{\lm_1}-a'_{2k+1}{\lm_2}\;
(=-h_{-2k-1})$.
  \end{enumerate}
 \label{lem2}
\end{lem1}
{\sl Proof.}\,\,
(i) First, we consider the case that $j$ is odd.
We will show $\sigma_{-j}(H_{-2k+1})=0$  
by the induction on $j$.
\vskip5pt
Case I) $j=1$.
\begin{align*}
 \sigma_{-1}(H_{-2k+1})& = h_{-1}-\lm_2\\ 
       & =\lm_2 - \lm_2\\
       & = 0.
\end{align*}
\vskip5pt
Case II) $j>1$.

\noindent
We assume  $\sigma_{-j+2}(H_{-2k+1})=0$.
Then, 
\begin{align*}
 \sigma_{-j}(H_{-2k+1})
        = & -\lan h_2,\lm\ran+
            \sum_{l\geq-j}\lan h_2,\al_{i_l}\ran x_l\\
        = & -\lan h_2,\lm\ran+
                           \sum_{l\geq-j+2}\lan h_2,\al_{i_l}\ran x_l
           +h_{-j} -c_2h_{-j+1}+h_{-j+2}\\
        = & \sigma_{-j+2}(H_{-2k+1})+h_{-j}+
          (-c_2a_{j-2}+a'_{j-3}){\lm_1}+
          (-c_2a_{j-1}+a'_{j-2}){\lm_2}\\
        = & h_{-j}-a'_{j-1}\lm_1-a'_{j}\lm_2 \q
           (\text{by  Fact \ref{af1}(ii),(iii)})\\
        = & 0.
\end{align*}
Next, we consider the case that $j$ is even. We will show
by the induction on $j$.
\vskip5pt
Case I) $j=2$.
\begin{align*}
 \sigma_{-2}(H_{-2k+1})& = h_{-2}-c_1h_{-1}-\lm_1\\ 
       & = h_{-2} -{\lm_1}-c_1{\lm_2}\\
       & =h_{-2} -a_1{\lm_1}-a_2{\lm_2} \\
       & = 0.
\end{align*}
\vskip5pt
Case II) $j>2$.

\noindent
We assume  $\sigma_{-j+2}(H_{-2k+1})=0$.
Then, 
\begin{align*}
 \sigma_{-j}(H_{-2k+1})
        = & -\lan h_1,\lm\ran+
            \sum_{l\geq-j}\lan h_1,\al_{i_l}\ran x_l\\
        = & -\lan h_1,\lm\ran+
                           \sum_{l\geq-j+2}\lan h_1,\al_{i_l}\ran x_l
           +h_{-j} -c_1h_{-j+1}+h_{-j+2}\\
        = & \sigma_{-j+2}(H_{-2k+1})+h_{-j}+
          (-c_1a'_{j-2}+a_{j-3}){\lm_1}+
          (-c_1a'_{j-1}+a_{j-2}){\lm_2}\\
        = & h_{-j}-a_{j-1}\lm_1-a_{j}\lm_2 \q
           (\text{by  Fact \ref{af1}(i),(iv)})\\
        = & 0.
\end{align*}

We can show (ii) by the similar calculation as (i).
% using Fact $\ref{af1}(ii),(iii)$. 
(iii), (iv) is obvious by (i), (ii).
\qed
%\hfill {\sl Q.E.D.}

 \begin{lem1}
If the condition of a weight $\lm$ is \,
 $\dfrac{a'_{2k-1}}{a'_{2k-2}}>\dfrac{\lm_1}{-\lm_2}
           \geq \dfrac{a_{2k}}{a_{2k-1}}$\,$(resp.\, \dfrac{a_{2k}}{a_{2k-1}}
>\dfrac{\lm_1}{-\lm_2}\geq \dfrac{a'_{2k+1}}{a'_{2k}}\,)$, then 
for any $j$ $($$1\leq j \leq k$$)$, 
 $H_{-2j+2}$, $H_{-2j+1}$\,$(resp.\,H_{-2j+1},\, H_{-2j}\,)$ 
are generated from $(\cd,0,0,t_{\lm},0,0,\cd)$.
\label{v}
\end{lem1}
{\sl Proof.}\,\,
We consider the former case $\dfrac{a'_{2k-1}}{a'_{2k-2}}
>\dfrac{\lm_1}{-\lm_2}\geq \dfrac{a_{2k}}{a_{2k-1}}$ since the latter
case is shown by the similar argument of the former case. 
By Remark $\ref{rem15}$, note that $h_{-1},\,h_{-2},\,\cd,h_{-2k+1}<
0$ and $h_{-2k}>0$.
We will show the lemma by the induction on $j$.
\vskip5pt
Case I) $j=1$.

\noindent
Let $v_0:=(\cd,0,0,t_{\lm},0,0,\cd)$. First, we show
$H_{-1}={\tilde{e}_2}^{-h_{-1}}v_0$. For any $s \geq 0$, $t \geq 1$, we have
$$ \sigma_{-1-2s}(v_0)=-\lm_2=-h_{-1},$$
$$\sigma_{2t}(v_0)=0.$$
By Lemma $\ref{lem1}$, we get 
$${\tilde{e}_2}^{-h_{-1}}v_0=(\cd,0,0,t_{\lm},h_{-1},0,\cd)=H_{-1}.$$
Note that ${\tilde{e}_2}^{-h_{-1}+1}v_0=0$. Next, we show
$H_{-2}={\tilde{e}_1}^{-h_{-2}}H_{-1}$. We will calculate the action of
${\tilde{e_1}}$ to $H_{-1}$. For any $s \geq 0$, $t \geq 1$, we have
$$ \sigma_{-2-2s}(H_{-1})=-c_1\lm_2-\lm_1=-h_{-2},$$
$$\sigma_{2t-1}(H_{-1})=0.$$
By Lemma $\ref{lem1}$, we get
$${\tilde{e}_1}^{-h_{-2}}H_{-1}=(\cd,0,0,t_{\lm},h_{-1},\,h_{-2},\cd)=H_{-2}.$$

Case II) $j>1$.

\noindent
We assume that $H_{-2j+3}$ is generated from 
 $(\cd,0,0,t_{\lm},0,0,\cd)$ 
by acting ${\tilde{e}_i}$'s.
First, we will calculate the action of $\tilde{e}_2$ to $H_{-2j+3}$.
By Lemma $\ref{lem2}(i)$, we have
$$0=\sigma_{-2j+3}(H_{-2j+3})=\sigma_{-2j+5}(H_{-2j+3})=
\cd=\sigma_{-1}(H_{-2j+3}).$$
By Lemma $\ref{lem1}(i),(ii)$, we have for any $s ,t \geq 1$,
$$\sigma_{-2j+3}(H_{-2j+3})>\sigma_{-2j+3-2s}(H_{-2j+3}),$$
$$\sigma_{2t}(H_{-2j+3})=0.$$
Therefore, we get 
$$\tilde{e}_2(H_{-2j+3})=0.$$

Next, we will calculate the action of $\tilde{e}_1$ to $H_{-2j+3}$. 
By Lemma $\ref{lem2}(iii)$, we have
$$\sigma_{-2j+2}(H_{-2j+3})=-a_{2j-3}\lm_1-a_{2j-2}\lm_2=-h_{-2j+2}
=:M>0$$
and
$$0=\sigma_{-2j+4}(H_{-2j+3})=\sigma_{-2j+6}(H_{-2j+3})
=\cd=\sigma_{-2}(H_{-2j+3}).$$
By Lemma $\ref{lem1}(i),(ii)$, we have for any $s ,t \geq 1$,
$$\sigma_{-2j+2}(H_{-2j+3})=\sigma_{-2j+2-2s}(H_{-2j+3}),$$
$$\sigma_{2t-1}(H_{-2j+3})=0.$$
Therefore, applying Lemma $\ref{lem1}(iii)$, we obtain
\begin{align*}
&{\tilde{e}_1}^{M}(H_{-2j+3})\;\;\,\,=H_{-2j+2},\\
&{\tilde{e}_1}^{M+1}(H_{-2j+3})={\tilde{e}_1}(H_{-2j+2})=0.
\end{align*}
%$${\tilde{e}_1}^{M}(H_{-2j+3})=H_{-2j+2}.$$
%$${\tilde{e}_1}^{M+1}(H_{-2j+3})={\tilde{e}_1}(H_{-2j+2})=0.$$

Finally, we will calculate the action of $\tilde{e}_2$ to $H_{-2j+2}$. 
By Lemma $\ref{lem2}(i)$, we have
$$\sigma_{-2j+1}(H_{-2j+2})=-a'_{2j-2}\lm_1-a'_{2j-1}\lm_2=-h_{-2j+1}
=:N>0$$
and
$$0=\sigma_{-2j+3}(H_{-2j+2})=\sigma_{-2j+4}(H_{-2j+2})
=\cd=\sigma_{-1}(H_{-2j+2}).$$
By Lemma $\ref{lem1}(i),(ii)$, we have for any $s ,t \geq 1$,
$$\sigma_{-2j+1}(H_{-2j+2})=\sigma_{-2j+1-2s}(H_{-2j+2}),$$
$$\sigma_{2t}(H_{-2j+2})=0.$$
Therefore, we obtain
$${\tilde{e}_2}^{N}(H_{-2j+2})=H_{-2j+1}.$$
( Note that ${\tilde{e}_2}^{N+1}(H_{-2j+2})={\tilde{e}_2}(H_{-2j+1})=0$.)
\qed
%\hfill {\sl Q.E.D.}
\vskip5pt

Now, we prove Theorem $\ref{hwv}(i)$ as below:

\noindent
In Lemma $\ref{v}$, we set $j=k$. Then
\begin{align*}
H_{-2k+1}=
(\cd,0,0,t_{\lm},\; h_{-1},\; h_{-2},\;,\cd,h_{-2k+2},\; h_{-2k+1},0,0,\cd)
\end{align*}
and
$${\tilde{e}_2}(H_{-2k+1})=0.$$
We will calculate the action of $\tilde{e}_1$ to $H_{-2k+1}$.
By Lemma $\ref{lem2}(i)$, we have 
$$0=\sigma_{-2k+2}(H_{-2k+1})=\sigma_{-2k+4}(H_{-2k+1})=
\cd=\sigma_{-2}(H_{-2k+1}).$$
In this case, by Lemma $\ref{lem2}(i)$ and Lemma $\ref{lem1}(ii)$, we have  
for any $s \geq 1$,
$$\sigma_{-2k}(H_{-2k+1})
=-a_{2k-1}\lm_1-a_{2k}\lm_2=\sigma_{-2k-2s}(H_{-2k+1}).$$
Since the $\lm$'s condition is
$\dfrac{a'_{2k-1}}{a'_{2k-2}}>\dfrac{\lm_1}{-\lm_2}
\geq \dfrac{a_{2k}}{a_{2k-1}}$, we have
 $-a_{2k-1}\lm_1-a_{2k}\lm_2 = -h_{-2k}\leq0$.
This shows
$$\tilde{e}_{1}(H_{-2k+1})=0$$
and $H_{-2k+1}$ is the highest weight vector of $B_0(\lm)$.

Similarly, if the $\lm$'s condition is 
$\dfrac{a_{2k}}{a_{2k-1}}>\dfrac{\lm_1}{-\lm_2}
           \geq \dfrac{a'_{2k+1}}{a'_{2k}}$, we can show that
 $H_{-2k}$ is the highest weight vector of $B_0(\lm)$.
\qed
%\hfill {\sl Q.E.D.}\\
\q\\

We can describe an explicit form of lowest weight vectors similarly to
describe highest weight vectors. For any $j,k \geq 1$, we define
\begin{align*}
&l_{2k-1}\,:=\;a_{2k-1}{\lm_1}\,+\,a_{2k-2}{\lm_2},\\
&l_{2k}\q\,:=\;a'_{2k}{\lm_1}\,+\,a'_{2k-1}{\lm_2},\\ 
&L_{j}\q\;:=
(\cd,0,0,\, l_j, \, l_{j-1}, \cd , l_2, \, l_1,t_{\lm},0,0,\cd).
\end{align*}
Now, we note that the following remark:
 \begin{rem1}
If the condition of a weight $\lm$ is $\;\dfrac{a_{2k}}{a_{2k-1}}\,<\,
\dfrac{\lm_1}{-\lm_2}\, \leq\, \dfrac{a'_{2k-1}}{a'_{2k}}$
(resp. $\dfrac{a'_{2k-1}}{a'_{2k}}\,<\,
\dfrac{\lm_1}{-\lm_2}\,\leq\, \dfrac{a_{2k}}{a_{2k+1}}$ ), then
$l_{1},\,l_{2},\,\cd,l_{2k-1}> 0$ and $l_{2k},\,l_{2k+1},\cd<0$ (resp.
 $l_{1},\,l_{2},\,\cd,l_{2k}> 0$ and $l_{2k+1},\,l_{2k+2},\cd<0$ ) by
Corollary $\ref{rem1}(ii)$.
\label{rem17} 
\end{rem1}
\begin{thm1}[Lowest weight vector of $B_0(\lm)$]
If the condition of a weight $\lm$ is $\dfrac{a_{2k}}{a_{2k-1}}\,<\,
\dfrac{\lm_1}{-\lm_2}\, \leq\, \dfrac{a'_{2k-1}}{a'_{2k}} $ 
$($resp. $\dfrac{a'_{2k-1}}{a'_{2k}}\,<\,
\dfrac{\lm_1}{-\lm_2}\,\leq\, \dfrac{a_{2k}}{a_{2k+1}}$$)$, 
then $L_{2k-1}$ $($resp. $L_{2k}$$)$ is the
   lowest weight vector of $B_0(\lm)$. 
\label{lwv}
 \end{thm1}
{\sl Proof.}\,\,This is shown by the similar way to the Theorem \ref{hwv}.
\qed
\vskip5pt
\noindent
By Corollary $\ref{rem1}$ and Theorem $\ref{hwv}$, we obtain the following
corollary.
\begin{cor1}
The existence conditions for the
highest $($or lowest $)$ weight vector of $B_0(\lm)$ is given by:
\begin{enumerate}
\item 
 If the $\lm$'s condition is\, $ \dfrac{\lm_1}{-{\lm_2}}>\alpha$, \,then 
$B_0(\lm)$  have the highest weight vector. 
\item
  If the $\lm$'s condition is\, $\beta >\dfrac{\lm_1}{-{\lm_2}}$, \,then  
$B_0(\lm)$ have the lowest weight vector. 
\end{enumerate}
\label{cor1}
\end{cor1}

 \begin{rem1}
Corollary $\ref{cor1}$ means that a weight $\lm$ is contained
in Tits cone.
\end{rem1}

\begin{rem1}
 \begin{enumerate}
  \item If $\ge$ is affine, then $\alpha = \beta$.
    In this case, if $\alpha = \beta = \dfrac{\lm_1}{-{\lm_2}}$, 
then $level(\lm)=0$
 $(\,level(\lm):=\langle c,\lm \rangle$, $c$ : center of $\ge\,)$.
  \item If $\ge$ is hyperbolic, then $\alpha > \beta$.
  In this case, if $\alpha \geq \dfrac{\lm_1}{-{\lm_2}} \geq \beta$, 
 $B_0(\lm)$ does not have the highest weight vector and lowest weight
vector.
 \end{enumerate}
\end{rem1}

%%%%%%%%% Section 5  %%%%%%%%%%%%%
\section{Polyhedral realization of $B_0(\lm)$ for affine and 
hyperbolic types of rank 2}
\setcounter{equation}{0}
\renewcommand{\theequation}{\thesection.\arabic{equation}}
%%%%%%%%%%%%%%%%%%%%

In this section, we shall describe the explicit form of $B_0(\lm)$ 
by polyhedral realization. We recall the definition of 
$\Xi_{\io}[\lm]$ and $\Sigma_{\io}[\lm]$:
 \begin{align*}
  \Xi_{\io}[\lm]:&=\{\bar{S}_{j_l}\cd \bar{S}_{j_1}(x_{j_0})\,
  :\,l\geq0,\,j_0,\cd,j_l\geq 1\}
  \\
  &\cup\{\bar{S}_{{-j}_k}\cd \bar{S}_{{-j}_1}(-x_{{-j}_0})\,
  :\,k\geq0,\,j_0,\cd,j_k \geq 1\},\\
 \Sigma_{\io}[\lm]
 :&=\{\vec x\in \ZZ^{\ify}_{\io}[\lm](\subset \QQ^{\ify}[\lm])\,:\,
\vp(\vec x)\geq 0\,\,{\rm for \,\,any }\,\,\vp\in \Xi_{\io}[\lm]\}.
\end{align*}

\noindent
First, we treat the case that $B_0(\lm)$ contains the highest 
weight vector.  We consider the following two cases of $\lm$'s condition:

(i) $\dfrac{a'_{2k-1}}{a'_{2k-2}}>\dfrac{\lm_1}{-\lm_2}
           \geq \dfrac{a_{2k}}{a_{2k-1}}$.
(ii) $\dfrac{a_{2k}}{a_{2k-1}}>\dfrac{\lm_1}{-\lm_2}
           \geq \dfrac{a'_{2k+1}}{a'_{2k}}$.
\vskip5pt
\noindent
We define
\begin{enumerate}
\item
If \;$\dfrac{a'_{2k-1}}{a'_{2k-2}}>\dfrac{\lm_1}{-\lm_2}
           \geq \dfrac{a_{2k}}{a_{2k-1}}$, 
\begin{align*}
  \Xi^{1}_{\io}[\lm]:=&
  \{\bar{S}_{{-j}_r}\cd \bar{S}_{{-j}_1}(x_{-m}-h_{-m})\,
  :\,r\geq1,\,j_r\geq1,\,1 \leq m \leq 2k-1\}
  \\
  \cup&\{\bar{S}_{{-j}_s}\cd \bar{S}_{{-j}_1}(x_{-2k})\,
  :\,s\geq1,\,j_s\geq1\},\\
\Sigma^{1}_{\io}[\lm]
 :=&\{\vec x\in \ZZ^{\ify}_{\io}[\lm]\,:\,
 \vp(\vec x)\geq 0\,\,{\rm for \,\,any }\,\,\vp\in
 \Xi^{1}_{\io}[\lm]\}.
\end{align*}
\item
If \;$\dfrac{a_{2k}}{a_{2k-1}}>\dfrac{\lm_1}{-\lm_2}
           \geq \dfrac{a'_{2k+1}}{a'_{2k}}$,
\begin{align*} 
  \Xi^{2}_{\io}[\lm]:=&
  \{\bar{S}_{{-j}_r}\cd
  \bar{S}_{{-j}_1}(x_{-m}-h_{-m})\,
  :\,r\geq1,\,j_r\geq1,\,1 \leq m \leq 2k\}
  \\
  \cup&\{\bar{S}_{{-j}_s}\cd \bar{S}_{{-j}_1}(x_{-2k-1})\,
  :\,s\geq1,\,j_s\geq1\},\\
\Sigma^{2}_{\io}[\lm]
 :=&\{\vec x\in \ZZ^{\ify}_{\io}[\lm]\,:\,
 \vp(\vec x)\geq 0\,\,{\rm for \,\,any }\,\,\vp\in
 \Xi^{2}_{\io}[\lm]\}.
\end{align*}
\end{enumerate}

\begin{thm1} 
 \label{main}
We fix an infinite sequence $\io = (\cd,2,1,t_{\lm},2,1,\cd)$.  
If the condition of a weight $\lm$ is 
\;$\dfrac{a'_{2k-1}}{a'_{2k-2}}>\dfrac{\lm_1}{-\lm_2}
           \geq \dfrac{a_{2k}}{a_{2k-1}}$ $(resp. 
\; \dfrac{a_{2k}}{a_{2k-1}}>\dfrac{\lm_1}{-\lm_2}
           \geq \dfrac{a'_{2k+1}}{a'_{2k}})$, 
$B_0(\lm)$ is described as below:
$$B_0(\lm) = \Sigma_{\io}[\lm]\cap\Sigma^{1}_{\io}[\lm]$$
$$(\text{resp. } B_0(\lm)= \Sigma_{\io}[\lm]\cap\Sigma^{2}_{\io}[\lm]).$$
\end{thm1}

\vskip5pt

{\sl Proof.}\,\,
We shall show this theorem in the former case (i) 
$\dfrac{a'_{2k-1}}{a'_{2k-2}}>\dfrac{\lm_1}{-\lm_2}
           \geq \dfrac{a_{2k}}{a_{2k-1}}$ since we can show in
the latter case (ii) $\dfrac{a_{2k}}{a_{2k-1}}>\dfrac{\lm_1}{-\lm_2}
           \geq \dfrac{a'_{2k+1}}{a'_{2k}}$
similarly to former case.
Now, $\Xi^{1}_{\io}[\lm]$ is closed by 
$\bar{S}_{k}$'s,  by Lemma \ref{mlem}
$\Sigma_{\io}[\lm]\cap\Sigma^{1}_{\io}[\lm]$ 
has a crystal structure unless it is empty. 
We will show that $\Sigma_{\io}[\lm]\cap\Sigma^{1}_{\io}[\lm]$ 
contains $\vec 0$ ,which implies that 
$\Sigma_{\io}[\lm]\cap\Sigma^{1}_{\io}[\lm]$ is non-empty,
and show that $\Sigma_{\io}[\lm]\cap\Sigma^{1}_{\io}[\lm]$ 
has the unique highest weight vector. 
First, we will show $\Sigma_{\io}[\lm]\cap\Sigma^{1}_{\io}[\lm]$ 
contains $\vec 0$.
For the purpose, we shall evaluate the constant term of
%By the similar to the proof of Theorem$\ref{thm1}$, we will show $\vec 0 \in
%\Sigma{'}_{\io}[\lm]$.
$\bar{S}_{-{j_k}}\cd\bar{S}_{-{j_1}}({x_{-k}})$. 

We set
\begin{equation*}
\vp^{(l)}_{-k}(x):=
\begin{cases}
\bar{S}_{-k+l-1}\cd\bar{S}_{-k+1}\bar{S}_{-k}(x_{-k})&(l\leq k),\\
\bar{S}_{-k+l}\cd\bar{S}_{1}\bar{S}_{-1}\bar{S}_{-2}
\cd\bar{S}_{-k+1}\bar{S}_{-k}(x_{-k})&(l>k).
\end{cases}
\end{equation*}
By the similar argument in $\cite{NZ}$ Lemma 4.2, 
we obtain the explicit form of
$\vp^{(l)}_{-k}(x)$ up to constant term as follows:

\begin{lem1}
Let $k$ be odd $($resp. even$)$, then 
\begin{align*}
%&&\Xi^{'}_{\io}[\lm]=\{\vp^{(l)}_{-k}(x)+h_{-k}\,|\,k \geq
%1,\,l\geq1\,\},\\
&\vp^{(l)}_{-k}(x)-\vp^{(l)}_{-k}(0)=
a'_{l+1}x_{l-k+\theta(l-k)}-a'_{l}x_{l-k+1+\theta(l-k+1)}\\
(resp.\;\;&\vp^{(l)}_{-k}(x)-\vp^{(l)}_{-k}(0)=
a_{l+1}x_{l-k+\theta(l-k)}-a_{l}x_{l-k+1+\theta(l-k+1)}\,)
\end{align*}
where
\[
\theta(x):=
\begin{cases}
1 &\text{if \,$x\geq0$,}\\
0 &\text{if \,$x<0$.}
\end{cases}
\]
\end{lem1} 

Now, we calculate the constant term in $\vp^{(l)}_{-k}$.
We consider the case that $k$ is odd.
For any $l \leq k-2$, we know that $\vp^{(l)}_{-k}$ has no costant
term by its definition. And we have
$$\vp^{(k-2)}_{-k}(x)=a'_{k-1}x_{-2}-a'_{k-2}x_{-1}.$$
By direct calculations, we obtain
\begin{eqnarray*}
\vp^{(k-1)}_{-k}&=&\bar{S}_{-2}(\vp^{(k-2)}_{-k})\\
                &=&a'_{k}x_{-1}-a'_{k-1}x_{1}+a'_{k-1}\lm_{1},\\
\vp^{(k)}_{-k}\q  &=&\bar{S}_{-1}\bar{S}_{-2}(\vp^{(k-2)}_{-k})\\
        &=&a'_{k+1}x_1-a'_{k}x_2+a'_{k-1}\lm_1+a'_{k}\lm_2. 
\end{eqnarray*}
For any $l> k$, since $\bar{S}_{-k+l}$ does not produce non-trivial
constant term, we have
$$\vp^{(k)}_{-k}(0)=\vp^{(l)}_{-k}(0).$$
Hence, we obtain the constant term of $\vp^{(l)}_{-k}(x)$:
\begin{equation*}
\vp^{(l)}_{-k}(0)=
\begin{cases} 
a'_{k-1}\lm_1+a'_{k}\lm_2 &\,\,(l\geq k),\\
a'_{k-1}\lm_1 &\,\,(l=k-1),\\
0 &\,\,(l\leq k-2).
\end{cases}
\end{equation*}
By remark \ref{rem15}, in the condition $(i)$ 
$\dfrac{a'_{2k-1}}{a'_{2k-2}}>\dfrac{\lm_1}{-\lm_2}
           \geq \dfrac{a_{2k}}{a_{2k-1}}$ 
we have $h_{-1},h_{-2},\cd,h_{-k} <0$ and $h_{-k-1}, h_{-k-2}, \cd >0$.
So, by the definition of $h_{-j}$, we have
$-h_{-j}+\vp^{(l)}_{-j}(0)\geq 0$
for any $j$ $(1 \leq j \leq k)$ and $\vp^{(l)}_{-j}(0)\geq 0$ for any 
$j$ $(j>k)$.
This shows that constant terms in all elements in $\Xi^{1}_{\io}[\lm]$
are non-negative and then $\vec 0$ is contained in 
$\Sigma_{\io}[\lm]\cap\Sigma^{1}_{\io}[\lm]$.

In the case that $k$ is even, we have
\begin{equation*}
\vp^{(l)}_{-k}(0)=
\begin{cases} 
a_{k-1}\lm_1+a_{k}\lm_2 &\,\,(l\geq k),\\
a_{k-1}\lm_1 &\,\,(l=k-1),\\
0 &\,\,(l\leq k-2)
\end{cases}
\end{equation*}
and $-h_{-j}+\vp^{(l)}_{-j}(0)\geq 0$ for any $j$ $(1 \leq j \leq
k-1)$
and $\vp^{(l)}_{-j}(0)\geq 0$ for any 
$j$ $(j>k)$.
 Therefore, 
$\vec 0$ is contained in 
$\Sigma_{\io}[\lm]\cap\Sigma^{1}_{\io}[\lm]$.
\vskip5pt

Next, we will show that $\Sigma_{\io}[\lm]\cap\Sigma^{1}_{\io}[\lm]$ 
has the unique highest weight vector $H_{-k}$.
Let $v_{\lm}:=(\cd,0,0,t_{\lm},x_{-1},x_{-2},\cd,x_{-k},\cd)$ 
be a highest weight vector in 
$\Sigma_{\io}[\lm]\cap\Sigma^{1}_{\io}[\lm]$, which satisfies:
\begin{equation}
 \begin{cases}
     x_{-j}- h_{-j} \geq 0  & (j \leq k),\\
       x_{-j} \geq 0     & (j >k).
 \end{cases}
\label{ineq}
\end{equation}
Note that the linear function $x_{-j}-h_{-j}$ $(j \leq k)$, 
$x_{-j}$ $(j>k)$ are generators 
of $\Xi^{1}_{\io}[\lm]$, and then
 any vector in
$\Sigma_{\io}[\lm]\cap\Sigma^{1}_{\io}[\lm]$
satisfies the inequality $(\ref{ineq})$.

We shall show that $v_{\lm}$ is uniquely determined and coincides with
$v_{\lm}$ by the induction on the index $j$ in the following two
cases: (I) $j \leq k$. (II) $j>k$.

Note that by the condition that $v_{\lm}$ is a highest weight vector, we have :
$$\sigma_{-j}(v_{\lm})\leq0 {\rm \,\,for \,\,any\,\, } j \geq 1.$$

CaseI) $j \leq k$.

\noindent
For $j=1$, we have 
$\sigma_{-1}(v_{\lm})=x_{-1}-\lm_2\leq0$.
By (\ref{ineq}), we also have 
$x_{-1}-h_{-1}=x_{-1}-\lm_2\geq0$,
which implies
\[
 x_{-1}=\lm_2=h_{-1}.
\]
Assume that for any $j'< j$, 
\begin{equation}
x_{-j'}=-h_{-j'}.
\label{ass2}
\end{equation}
Now, by Lemma $\ref{lem2}$, note that we know that 
\begin{eqnarray*}
  x_0 =(\cd,0,0,t_{\lm},h_{-1},
    h_{-2},\cd,h_{-2k+1},\cd)\,\,
\end{eqnarray*}
is one of the highest weight vectors of 
$\ZZ^{\ify}_{\io}[\lm]$ and satisfies the Lemma \ref{lem2}.

Let us determine $x_{-j}$.
By the assumption (\ref{ass2}) and Lemma $\ref{lem2}(iii)$, we have
\begin{equation}
\sigma_{-j}(v_{\lm})=x_{-j}-h_{-j} \leq 0.
\label{xd2}
\end{equation}
By $(\ref{ineq})$, we have 
\[
 x_{-j}-h_{-j} \geq 0.
\]
Therefore, we obtain $x_{-j}=h_{-j}$.
\vskip5pt
CaseII) $j>k$.

\noindent
We note that
\begin{equation}
\sigma_{-1}(v_{\lm})= \sigma_{-2}(v_{\lm})= \cd = \sigma_{-k}(v_{\lm})=0
\end{equation}
by the assumption (\ref{ass2}) and Lemma $\ref{lem2}(iii)$.
For $j=k+1$, 
by Remark \ref{rem15},
 we have $h_{-k-1}>0$. 
In this case, since $x_{-k-1} \leq 0$, we know that
$\sigma_{-k-1}(v_{\lm})$ is non-positive automatically as follows:
\[
\sigma_{-k-1}(v_{\lm})=x_{-k-1}-h_{-k-1}\leq 0.
\]
But, $x_{-k-1}$ is a generator of $\Xi^{1}[\lm]$ and then
$$x_{-k-1} \geq 0.$$
This shows $x_{-k-1}=0$.

Now, assume that for any $j'<j$ $(k<j')$,
$$x_{-j'}=0.$$
Under the assumption, we have
\begin{equation*}
\begin{cases}
\sigma_{-1}(v_{\lm})= \sigma_{-2}(v_{\lm})= \cd =
\sigma_{-k}(v_{\lm})=0,\\
\sigma_{-j}(v_{\lm})=\sigma_{-k-1}(v_{\lm})=-h_{-k-1}\leq 0 
&(j : \text{even}),\\
\sigma_{-j}(v_{\lm})=0 &(j : \text{odd}).
\end{cases}
\end{equation*}
In this case, $x_{-j}$ is a generator of $\Xi^{1}[\lm]$ and then
$$x_{-j} \geq 0.$$
This shows $x_{-j}=0$.

Therefore, we obtain 
\begin{equation}
  x_{-j}=\begin{cases}
      h_{-j} & (j \leq k),\\
       0     & (j >k).
 \end{cases}
\end{equation}
Now, we know that $v_{\lm}=H_{-k}$ is the unique highest weight 
vector in $\Sigma_\io[\lm]\cap\Sigma^{1}_{\io}[\lm]$. 
Since $B_0(\lm)$ contains 
the unique highest weight vector \cite{K4}, 
$H_{-k}$ must be the unique highest weight vector in $B_0(\lm)$.
\qed

\vskip5pt
We show the results of the explicit form of $\Xi$ as follows: 
\begin{align*} 
\Xi_{\io}[\lm]\,&=\,
\left\{
\begin{array}{lr}
{a_m}x_m - a_{m-1}x_{m+1} & (m\geq1)\\
-a'_{m}x_{-m}+a'_{m-1}x_{-m-1} & (m\geq1)
\end{array}
\right\},\\
\q\\
\Xi^{1}_{\io}[\lm]&=
\left\{
\begin{array}{lr}
a'_{i+1}x_{-2j+i+1}-a'_{i}x_{-2j+i+2}-a'_{2j-2}\lm_1
-a'_{2j-1}\lm_2
&(1\leq j \leq k,\; 0\leq i \leq 2j-3)\\
a'_{2j-1}x_{-1}-a'_{2j-2}x_{1}-a'_{2j-1}\lm_2
&(1\leq j \leq k,\;i=2j-2)\\
a'_{i+1}x_{-2j+i+2}-a'_{i}x_{-2j+i+3}
&(1\leq j \leq k,\;i\geq 2j-1)\\
a_{i+1}x_{-2j+i+2}-a_{i}x_{-2j+i+3}-a_{2j-3}\lm_1-a_{2j-2}\lm_2
&(1\leq j \leq k,\;0\leq i \leq 2j-4)\\
a'_{i+1}x_{-2j+i+2}-a'_{i}x_{-2j+i+3}+a'_{2j-2}\lm_1
+a'_{2j-1}\lm_2 
&(j\geq k+1,\;i \geq 2j-1)\\
a_{i+1}x_{-2j+i+2}-a_{i}x_{-2j+i+3}
&(j\geq k+1,\;0 \leq i \leq 2j-4)\\
a_{2j-2}x_{-1}-a_{2j-3}x_{1}+a_{2j-3}\lm_1 
&(j\geq k+1,\;i=2j-3)\\
a_{i+1}x_{-2j+i+3}-a_{i}x_{-2j+i+4}+a_{2j-3}\lm_1+a_{2j-2}\lm_2 
&(j\geq k+1,\;i \geq 2j-2)
\end{array}
\right\},
\end{align*}
\begin{align*}
\Xi^{2}_{\io}[\lm]&=
\left\{
\begin{array}{lr}
a_{i+1}x_{-2j+i}-a_{i}x_{-2j+i+1}-a_{2j-1}\lm_1
-a_{2j}\lm_2
&(1\leq j \leq k,\; 0\leq i \leq 2j-2)\\
a_{2j}x_{-1}-a_{2j-1}x_{1}-a_{2j}\lm_2
&(1\leq j \leq k,\;i=2j-1)\\
a_{i+1}x_{-2j+i+1}-a_{i}x_{-2j+i+2}
&(1\leq j \leq k,\;i\geq 2j)\\
a'_{i+1}x_{-2j+i+1}-a'_{i}x_{-2j+i+2}-a'_{2j-2}\lm_1
-a'_{2j-1}\lm_2
&(1\leq j \leq k,\;0\leq i \leq 2j-3)\\
a_{i+1}x_{-2j+i+3}-a_{i}x_{-2j+i+4}+a_{2j-3}\lm_1
+a_{2j-2}\lm_2 
&(j\geq k+1,\;i \geq 2j-2)\\
a'_{i+1}x_{-2j+i+3}-a'_{i}x_{-2j+i+4}
&(j\geq k+1,\;0 \leq i \leq 2j-5)\\
a'_{2j-3}x_{-1}-a'_{2j-4}x_{1}+a'_{2j-4}\lm_1 
&(j\geq k+1,\;i=2j-4)\\
a'_{i+1}x_{-2j+i+4}-a'_{i}x_{-2j+i+5}+a'_{2j-4}\lm_1
+a'_{2j-3}\lm_2 
&(j\geq k+1,\;i \geq 2j-3)
\end{array}
\right\}.
\end{align*}
Note that constant terms of all linear forms in 
$\Xi_{\io}[\lm]$, $\Xi^{1}_{\io}[\lm]$ and $\Xi^{2}_{\io}[\lm]$
are nonnegative.

Next, we treat the case that $B_0(\lm)$ contains the lowest 
weight vector.  We consider the following two cases of $\lm$'s condition:

(i) $\dfrac{a'_{2k-1}}{a'_{2k}}\geq\dfrac{\lm_1}{-\lm_2}
           > \dfrac{a_{2k}}{a_{2k-1}}$.
(ii) $\dfrac{a_{2k}}{a_{2k+1}}\geq\dfrac{\lm_1}{-\lm_2}
           > \dfrac{a'_{2k-1}}{a'_{2k}}$.

\vskip5pt
\noindent 
We define
\begin{enumerate}
\item
If \;$\dfrac{a'_{2k-1}}{a'_{2k}}\geq\dfrac{\lm_1}{-\lm_2}
           > \dfrac{a_{2k}}{a_{2k-1}}$,
\begin{align*} 
  \Xi^{3}_{\io}[\lm]:=&
  \{\bar{S}_{{j}_r}\cd
  \bar{S}_{{j}_1}(-x_m+l_m)\,
  :\,r\geq1,\,j_r\geq1,\,1 \leq m \leq 2k-1\}
  \\
  \cup&\{\bar{S}_{{j}_s}\cd \bar{S}_{{j}_1}(-x_{2k})\,
  :\,s\geq1,\,j_s\geq1\},\\
\Sigma^{3}_{\io}[\lm]
 :=&\{\vec x\in \ZZ^{\ify}_{\io}[\lm]\,:\,
 \vp(\vec x)\geq 0\,\,{\rm for \,\,any }\,\,\vp\in
 \Xi^{3}_{\io}[\lm]\}.
\end{align*}
\item
If \;$\dfrac{a_{2k}}{a_{2k+1}}\geq\dfrac{\lm_1}{-\lm_2}
           > \dfrac{a'_{2k-1}}{a'_{2k}}$,
\begin{align*}
  \Xi^{4}_{\io}[\lm]:=&
  \{\bar{S}_{{j}_r}\cd
  \bar{S}_{{j}_1}(-x_m+l_m)\,
  :\,r\geq1,\,j_r\geq1,\,1 \leq m \leq 2k\}
  \\
  \cup&\{\bar{S}_{{j}_s}\cd \bar{S}_{{j}_1}(-x_{2k+1})\,
  :\,s\geq1,\,j_s\geq1\},\\
\Sigma^{4}_{\io}[\lm]
 :=&\{\vec x\in \ZZ^{\ify}_{\io}[\lm]\,:\,
 \vp(\vec x)\geq 0\,\,{\rm for \,\,any }\,\,\vp\in
 \Xi^{4}_{\io}[\lm]\}.
\end{align*}
\end{enumerate}

\begin{thm1} 
 \label{main2}
We fix an infinite sequence $\io = (\cd,2,1,t_{\lm},2,1,\cd)$.  
If the condition of a weight $\lm$ is\; 
$\dfrac{a'_{2k-1}}{a'_{2k}}\geq\dfrac{\lm_1}{-\lm_2}
           > \dfrac{a_{2k}}{a_{2k-1}}$ $(resp. 
\; \dfrac{a_{2k}}{a_{2k+1}}\geq\dfrac{\lm_1}{-\lm_2}
           > \dfrac{a'_{2k-1}}{a'_{2k}})$, 
$B_0(\lm)$ is described as below:
$$B_0(\lm) = \Sigma_{\io}[\lm]\cap\Sigma^{3}_{\io}[\lm]$$
$$(\text{resp. } B_0(\lm)= \Sigma_{\io}[\lm]\cap\Sigma^{4}_{\io}[\lm]).$$
\end{thm1}

\vskip5pt

{\sl Proof.}\,\,This is the same as Theorem \ref{main}.
\qed
\vskip5pt

We show the results of the explicit form of $\Xi$ as follows: 
\begin{align*} 
\Xi_{\io}[\lm]\,&=\,
\left\{
\begin{array}{lr}
{a_m}x_m - a_{m-1}x_{m+1} & (m\geq1)\\
-a'_{m}x_{-m}+a'_{m-1}x_{-m-1} & (m\geq1)
\end{array}
\right\},\\
\q\\
\Xi^{3}_{\io}[\lm]&=
\left\{
\begin{array}{lr}
-a_{i+1}x_{2j-i-1}+a_{i}x_{2j-i-2}+a_{2j-1}\lm_1
+a_{2j-2}\lm_2
&(1\leq j \leq k,\; 0\leq i \leq 2j-3)\\
-a_{2j-1}x_{1}+a_{2j-2}x_{-1}+a_{2j-1}\lm_1
&(1\leq j \leq k,\;i=2j-2)\\
-a_{i+1}x_{2j-i-2}+a_{i}x_{2j-i-3}
&(1\leq j \leq k,\;i\geq 2j-1)\\
-a'_{i+1}x_{2j-i-2}+a'_{i}x_{2j-i-3}+a'_{2j-2}\lm_1+a'_{2j-3}\lm_2
&(1\leq j \leq k,\;0\leq i \leq 2j-4)\\
-a_{i+1}x_{2j-i-2}+a_{i}x_{2j-i-3}-a_{2j-1}\lm_1
-a_{2j-2}\lm_2 
&(j\geq k+1,\;i \geq 2j-1)\\
-a'_{i+1}x_{2j-i-2}+a'_{i}x_{2j-i-3}
&(j\geq k+1,\;0 \leq i \leq 2j-4)\\
-a'_{2j-2}x_{1}+a'_{2j-3}x_{-1}-a'_{2j-3}\lm_2 
&(j\geq k+1,\;i=2j-3)\\
-a'_{i+1}x_{2j-i-3}+a'_{i}x_{2j-i-4}-a'_{2j-2}\lm_1
-a'_{2j-3}\lm_2 
&(j\geq k+1,\;i \geq 2j-2)
\end{array}
\right\},
\end{align*}
\begin{align*}
\Xi^{4}_{\io}[\lm]&=
\left\{
\begin{array}{lr}
-a'_{i+1}x_{2j-i}+a'_{i}x_{2j-i-1}+a'_{2j}\lm_1
+a'_{2j-1}\lm_2
&(1\leq j \leq k,\; 0\leq i \leq 2j-2)\\
-a'_{2j}x_{1}+a'_{2j-1}x_{-1}+a'_{2j}\lm_1
&(1\leq j \leq k,\;i=2j-1)\\
-a'_{i+1}x_{2j-i-1}+a'_{i}x_{2j-i-2}
&(1\leq j \leq k,\;i\geq 2j)\\
-a_{i+1}x_{2j-i-1}+a_{i}x_{2j-i-2}+a_{2j-1}\lm_1+a_{2j-2}\lm_2
&(1\leq j \leq k,\;0\leq i \leq 2j-3)\\
-a'_{i+1}x_{2j-i-1}+a'_{i}x_{2j-i-2}-a'_{2j}\lm_1
-a'_{2j-1}\lm_2 
&(j\geq k+1,\;i \geq 2j)\\
-a_{i+1}x_{2j-i-1}+a_{i}x_{2j-i-2}
&(j\geq k+1,\;0 \leq i \leq 2j-3)\\
-a_{2j-1}x_{1}+a_{2j-2}x_{-1}-a_{2j-2}\lm_2 
&(j\geq k+1,\;i=2j-2)\\
-a_{i+1}x_{2j-i-2}+a_{i}x_{2j-i-3}-a_{2j-1}\lm_1
-a_{2j-2}\lm_2 
&(j\geq k+1,\;i \geq 2j-1)
\end{array}
\right\}.
\end{align*}

%%%%%%%%% Section 6   %%%%%%%%%%
\appendix
\section{Highest and lowest weight vector of $B_0(\lm)$ 
for classical types 
of rank 2}

\setcounter{equation}{0}
\renewcommand{\theequation}{\thesection.\arabic{equation}}
%%%%%%%%%%%%%%%%%%%%

In this appendix, we consider the case of classical types $A_2$, $B_2$
%, $C_2$ 
and $G_2$. Deference of affine or hyperbolic type is that the 
integer $a_l$ may
not positive. We need to treat classical types 
case by case. If we obtain the
explicit form of the highest (or lowest) weight vector of $B_0(\lm)$ 
in $\ZZ^{\ify}_{\io}[\lm]$, we can describe the explicit form of 
$B_0(\lm)$ by the same method of previous section. 
In particular, the form of $\Xi$ is the same.
So, we describe the explicit form of the highest (or lowest) weight vector
of $B_0(\lm)$ in the following subsection. 

%%%% 7.1 %%%%%
\subsection{$A_2$ case}

In this case, we have $c_1=c_2$ and then $a_l=a'_l$. 
We calculate integer $a_l$ till $a_l$ is zero
(c.f. Example \ref{exa}). Then,
\begin{eqnarray*}
a_1=1,\,a_2=1,\,a_3=0.
\end{eqnarray*}
Therefore, the sequence in Corollary \ref{rem1} is
$$1>0.$$
If $\dfrac{\lm_1}{-\lm_2} > 1$, then $h_1< 0$, $h_2 > 0$ and
$l_1$, $l_2 > 0$.
If $\dfrac{\lm_1}{-\lm_2} = 1$, then $h_1< 0$, $h_2= 0$ and
$l_1>0$, $l_2= 0$.
If $1 > \dfrac{\lm_1}{-\lm_2} > 0$, then $h_1$, $h_2 < 0$
and $l_1>0$, $l_2 <0$.

Consequently, we describe the highest (or lowest) weight vector of 
$B_0(\lm)$ as follows:
  \begin{enumerate}
  \item If \,$\dfrac{\lm_1}{-\lm_2} > 1$, 
then $H_{-1}$ is the
   highest weight vector and $L_2$ is the lowest weight vector
   of $B_0(\lm)$.
  \item If \,$\dfrac{\lm_1}{-\lm_2} = 1$, 
then $H_{-1}$ is the
   highest weight vector and $L_1$ is the lowest weight vector
   of $B_0(\lm)$.
  \item If \,$1 > \dfrac{\lm_1}{-\lm_2} > 0$, 
then $H_{-2}$ is the
   highest weight vector and $L_1$ is the lowest weight vector 
of $B_0(\lm)$.
  \end{enumerate}

%%%% 7.2 %%%%%
\subsection{$B_2$ case}

In this case, we have $c_1=2,c_2=1$ and then
\begin{eqnarray*}
a_1=1,\,a_2=2,\,a_3=1,\,a_4=0,\\
a'_1=1,\,a'_2=1,\,a'_3=1,\,a'_4=0.
\end{eqnarray*}
The sequence in Corollary \ref{rem1} is
$$2>1>0.$$
Therefore, we obtain a following table 1($H$: highest weight vector, 
$L$: Lowest weight vector):

\begin{table}[h]
\caption{}
\begin{center}
\renewcommand{\arraystretch}{1.5}
\begin{tabular}{@{}l@{}}
\noalign{\hrule height0.4pt}
\begin{tabular}{|c|c|c|}
$\lm$'s condition & $H$ & $L$ \\
\hline
$\frac{\lm_1}{-\lm_2} > 2$ & $H_{-1}$ & $L_3$\q \\
$\frac{\lm_1}{-\lm_2} = 2$ & $H_{-1}$ & $L_2$\q \\
$2 > \frac{\lm_1}{-\lm_2} > 1$ & $H_{-2}$ & $L_2$\q\\
$ \frac{\lm_1}{-\lm_2} = 1$ & $H_{-2}$ & $L_1$\q\\
$1 > \frac{\lm_1}{-\lm_2} > 0$ & $H_{-3}$ & $L_1$\q
\end{tabular}\\
\noalign{\hrule height0.4pt}
\end{tabular}
\end{center}
\end{table}

\subsection{$G_2$ case}

In this case, we have $c_1=1,c_2=3$ and then
\begin{eqnarray*}
a_1=1,\,a_2=1,\,a_3=2,\,a_4=1,\,a_5=1,\,a_6=0,\\
a'_1=1,\,a'_2=3,\,a'_3=2,\,a'_4=3,\,a'_5=1,\,a'_6=0.
\end{eqnarray*}
The sequence in Corollary \ref{rem1} is
$$1>\frac{2}{3}>\frac{1}{2}>\frac{1}{3}>0.$$
Therefore, we obtain a following table 2:
\begin{table}[h]
\caption{}
\begin{center}
\renewcommand{\arraystretch}{1.5}
\begin{tabular}{@{}l@{}}
\noalign{\hrule height0.4pt}
\begin{tabular}{|c|c|c|}
$\lm$'s condition & $H$ & $L$ \\
\hline
$\frac{\lm_1}{-\lm_2} > 1$ & $H_{-1}$ & $L_5$\q \\
$\frac{\lm_1}{-\lm_2} = 1$ & $H_{-1}$ & $L_4$\q \\
$1 > \frac{\lm_1}{-\lm_2} > \frac{2}{3}$ & $H_{-2}$ & $L_4$\q\\
$ \frac{\lm_1}{-\lm_2} = \frac{2}{3}$ & $H_{-2}$ & $L_3$\q\\
$\frac{2}{3} > \frac{\lm_1}{-\lm_2} > \frac{1}{2}$ & $H_{-3}$ &
$L_3$\q\\
$\frac{\lm_1}{-\lm_2} = \frac{1}{2}$ & $H_{-3}$ & $L_2$\q\\
$\frac{1}{2} > \frac{\lm_1}{-\lm_2} > \frac{1}{3}$ & $H_{-4}$ &
$L_2$\q\\
$ \frac{\lm_1}{-\lm_2} = \frac{1}{3}$ & $H_{-4}$ & $L_1$\q\\
$\frac{1}{3} > \frac{\lm_1}{-\lm_2} > 0$ & $H_{-5}$ & $L_1$
\end{tabular}\\
\noalign{\hrule height0.4pt}
\end{tabular}
\end{center}
\end{table}
%  \begin{enumerate}
%  \item If \,$\frac{\lm_1}{-\lm_2} > 1$, 
%then $H_{-1}$ is the
%   highest weight vector and $L_5$ is the lowest weight vector
%   of $B_0(\lm)$.
%  \item If \,$\frac{\lm_1}{-\lm_2} = 1$, 
%then $H_{-1}$ is the
%   highest weight vector and $L_4$ is the lowest weight vector
%   of $B_0(\lm)$.
%  \item If \,$1 > \frac{\lm_1}{-\lm_2} > \frac{2}{3}$, 
%then $H_{-2}$ is the
%   highest weight vector and $L_4$ is the lowest weight vector 
%of $B_0(\lm)$,
%\item If \,$ \frac{\lm_1}{-\lm_2} = \frac{2}{3}$, 
%then $H_{-2}$ is the
%   highest weight vector and $L_3$ is the lowest weight vector 
%of $B_0(\lm)$,
%\item If \,$\frac{2}{3} > \frac{\lm_1}{-\lm_2} > \frac{1}{2}$, 
%then $H_{-3}$ is the
%   highest weight vector and $L_3$ is the lowest weight vector 
%of $B_0(\lm)$,
%\item If \,$ \frac{\lm_1}{-\lm_2} = \frac{1}{2}$, 
%then $H_{-3}$ is the
%   highest weight vector and $L_2$ is the lowest weight vector 
%of $B_0(\lm)$,
%\item If \,$\frac{1}{2} > \frac{\lm_1}{-\lm_2} > \frac{1}{3}$, 
%then $H_{-4}$ is the
%   highest weight vector and $L_2$ is the lowest weight vector 
%of $B_0(\lm)$,
%\item If \,$ \frac{\lm_1}{-\lm_2} = \frac{1}{3}$, 
%then $H_{-4}$ is the
%   highest weight vector and $L_1$ is the lowest weight vector 
%of $B_0(\lm)$,
%\item If \,$\frac{1}{3} > \frac{\lm_1}{-\lm_2} > 0$, 
%then $H_{-5}$ is the
%   highest weight vector and $L_1$ is the lowest weight vector 
%of $B_0(\lm)$.
%  \end{enumerate}

\noindent
{\large\bf Acknowledgment}

The author is grateful to Professor T.Nakashima for his
support and encouragement during the course of this work. He is also
thankful to T.Yokonuma, K.Shinoda, K.Gomi, Y.Koga and H.Miyachi for helpful
comments and kindly advices. 
%He also thanks K.Handa, M.Katsuta,
%K.Suzuki and K.Yasuda for their various opinions.

%%%%%%%%%%%  reference  %%%%%%%%%%%%%%

\end{document}